\documentclass[12pt]{article}
\usepackage{amsmath, amsthm, amssymb}
\usepackage{mathrsfs}
\usepackage{bm}
\usepackage{natbib}
\usepackage[dvips]{graphicx}
\numberwithin{equation}{subsection}
\newtheorem{thm}{Theorem}
\newtheorem{lem}{Lemma}[subsection]

\newtheorem{cor}{Corollary}
\newtheorem{definition}{Definition}
\newtheorem{remark}{Remark}
\newtheorem{example}{Example}

\newtheorem{exercise}{Exercise}
\newtheorem{report}{Report}
\newtheorem{assumption}{Assumption}

\newcommand{\Proof}{\noindent {\bf Proof: }\ }
\newcommand{\bin}{\mathrm{Bin}}

\newcommand{\prob}{\mathrm{Prob}}
\newcommand{\abs}[1]{\vert #1 \vert}
\newcommand{\abslr}[1]{\left\vert #1 \right\vert}

\newcommand{\real}{\mathbb{R}}
\newcommand{\naturalnum}{\mathbb{N}}
\newcommand{\dist}{\mathrm{dist}}
\newcommand{\indicator}{\mathrm{1}}
\newcommand{\rmd}{\mathrm{d}}
\newcommand{\rmor}{\;\mathrm{or}\;}
\newcommand{\io}{{i.o.}}

\newcommand{\vr}[1]{\boldsymbol{#1}}

\newcommand{\argsup}{\mathrm{argsup}}

\newcommand{\scra}{\mathscr{A}}
\newcommand{\scrb}{\mathscr{B}}
\newcommand{\scrf}{\mathscr{G}}

\newcommand{\scrk}{\mathscr{K}}

\newcommand{\scrm}{\mathscr{K}}
\newcommand{\bnd}{\tilde{d}}
\newcommand{\cnd}{d}

\title{Strong consistency of the maximum likelihood estimator for finite mixtures of location-scale distributions 
when penalty is imposed on the ratios of the scale parameters}
\author{Kentaro Tanaka\footnote{Department of Industrial Engineering and Management, Tokyo Institute of Technology, 
2-12-1 O-okayama, Meguro-ku, Tokyo 152-8552 JAPAN, E-mail:\,tanaka.k.al@m.titech.ac.jp}}

\begin{document}

\maketitle

\begin{abstract}
 In finite mixtures of location-scale distributions, 
 if there is no constraint or penalty on the parameters, 
 then the maximum likelihood estimator does not exist because the likelihood is unbounded. 
 To avoid this problem, we consider a penalized likelihood, 
 where the penalty is a function of the minimum of the ratios of the scale parameters and the sample size. 
 It is shown that the penalized maximum likelihood estimator is strongly consistent. 
 We also analyze the consistency of a penalized maximum likelihood estimator 
 where the penalty is imposed on the scale parameters themselves. 
\end{abstract}

\section{Introduction}

In this paper, we prove the strong consistency of a penalized maximum likelihood estimate 
for finite mixtures of univariate location-scale distributions 
generalizing the results in \citet{CRI2003}. 
As a special case of this result, we solve an open problem posed by \citet{H1985}. 

As stated in \citet{D1969}, 
because the likelihood function for finite mixtures of location-scale distributions is unbounded, 
the maximum likelihood estimator does not exist. 
To see that, we consider a simple case that the model consists of mixtures of two normal distributions 
$\alpha_{1}\phi(x;\mu_{1}, \sigma_{1}) + \alpha_{2}\phi(x;\mu_{2}, \sigma_{2})$ 
and assume that we obtain an i.i.d. \hspace{-.5mm}sample $X_{1}, X_{2}, \dots, X_{n}$ from the true distribution. 
If we set $\mu_{1} = X_{1}$ and $\sigma_{1} \rightarrow 0$, 
then the likelihood tends to infinity as $\sigma_{1}$ goes to zero. 
Hence the likelihood function is unbounded. 

A straightforward approach to this problem is to bound the minimum of the variances of the components from below by a positive constant. 
By using theorem 6 in \citet{R1981}, 
we can show that the maximum likelihood estimator restricted to a compact subset of the parameter space is strongly consistent 
if the subset contains the true parameter. 

Another approach is penalized maximum likelihood estimation. 
However, if the penalty is not appropriate, then the likelihood function is unbounded. 
\citet{CRI2003} considered the case that the penalties are imposed on the variances themselves 
and proved the consistency of the penalized maximum likelihood estimator. 
The results given in \citet{CRI2003} are very useful for estimating the parameters of mixture of normal densities  
because the assumptions for the penalty are easy to check 
and the implementation of their method is also easy. 
In this paper, we extend their consistency result to the case that 
the components of mixtures are not normal densities and the penalty depends on the sample size $n$. 

In normal mixture distributions, 
\citet{H1985} considered the following constraints to avoid the divergence of the likelihood. 
\begin{equation}
 \min_{m, m'}\frac{\sigma_{m}}{\sigma_{m'}} \geq b
  \label{eq:HathawayConst}
\end{equation}
This bounds the minimum of the ratios of the variances of the components by a constant. 
He showed that the strong consistency of the maximum likelihood estimator holds 
if the true distribution satisfies the constraint represented by equation (\ref{eq:HathawayConst}). 
Intuitively, a stronger constraint must be enforced for a smaller sample size to avoid the divergence of the likelihood, 
because a component with a very small variance can only have a large contribution to a single observation at most. 
Therefore, it seems that the constraint under which the consistency holds can be weakened as the sample size increases. 
This intuition leads to the following two questions: 
\begin{itemize}
 \item Is it possible to let $b$ decrease to zero as the sample size $n$ increases to infinity while maintaining consistency?
 \item If it is possible, then at what rate can $b$ be decreased to zero?
\end{itemize}
These questions are mentioned in \citet{H1985}, \citet{M2000}, and treated as unsolved problems. 

This topic is closely related to a sieve method. (See \citet{G1981} and \citet{G1982}. ) 
For normal mixture distributions, 
the convergence rate of the maximum likelihood estimator based on sieve method is studied in \citet{GW2000} and \citet{GV2001}. 
In \citet{TT2003-35}, for mixtures of location-scale distributions, 
we showed the strong consistency of the maximum likelihood estimator 
in the case that the scale parameters themselves are bounded from below by 
$c_{n} = e^{-n^{d}} , (0<d<1)$. 
However, we could not solve the original questions when 
constraints are imposed on the minimum of the ratios of the variances of the components. 

In this paper, we solve the questions treated above in a more general and unified framework. 
For mixtures of location-scale distributions, 
we consider a penalized likelihood, where the penalty is a function of the minimum of the ratios of the scale parameters 
and the sample size $n$. 
The effect of the penalty becomes stronger as the minimum of the ratios of the scale parameters decreases to zero. 
Note that the penalty can depend on the sample size $n$. 
We can weaken the effect of the penalty as the sample size $n$ increases to infinity. 
In Theorem~\ref{thm:penalt-ratio-scale}, we show that the consistency holds for the penalized maximum likelihood estimator. 
In Corollary~\ref{cor:const-ratio-scale}, the solutions to the questions mentioned in \citet{H1985}, \citet{M2000} are obtained 
as special cases of Theorem~\ref{thm:penalt-ratio-scale}. 
We also analyze the consistency of a penalized maximum likelihood estimator 
in which the penalties are imposed on the scale parameters themselves. 
The result obtained in Theorem~\ref{thm:penalt-scale-param} is a generalization of Corollary 1 of \citet{CRI2003}. 

Throughout this paper, we assume that the true distribution is a mixture of location-scale distributions 
and the number of components of the true distribution is known. 

The organization of this paper is as follows. 
Section~\ref{sec:preliminaries} describes notation and regularity conditions. 
The main results are stated in section~\ref{sec:main-results}. 
Section~\ref{sec:proofs} is devoted to the proofs. 
We end this paper by concluding remarks in section~\ref{sec:conclusion}. 

\section{Preliminaries}
\label{sec:preliminaries}
\subsection{Notation}

Mixture of $M$ location-scale densities are written in the form 
\begin{equation}
 f(x;\theta) \equiv \sum_{m=1}^{M} \alpha_{m}f_{m}(x;\mu_{m}, \sigma_{m}). 
 \nonumber 
\end{equation}
The mixing weights $\alpha_{1},\dots,\alpha_{M}$ 
have to satisfy $\alpha_{m} \geq 0 ,\, \sum_{m=1}^{M}\alpha_{m}=1$. 
We assume that the components 
$f_{1}(x;\mu_{1}, \sigma_{1}),\dots,f_{M}(x;\mu_{M}, \sigma_{M})$ 
are location-scale densities i.e. they satisfy 
\begin{equation}
 f_{m}(x;\mu_{m}, \sigma_{m}) = \frac{1}{\sigma_{m}}f_{m}(\frac{x-\mu_{m}}{\sigma_{m}}; 0, 1) 
  \quad , \quad 
  1 \leq m \leq M, 
  \nonumber 
\end{equation}
where $\mu_{m}$ and $\sigma_{m}$ are location parameters and scale parameters respectively. 
We abbreviate $(\alpha_{1},\mu_{1},\sigma_{1},\dots,\alpha_{M},\mu_{M},\sigma_{M})$ as $\theta$, 
and $(\mu_{m},\sigma_{m})$ as $\theta_{m}$. 
We denote the true parameter by $\theta_{0}$. 

Let $\Omega_{m} \equiv \{(\mu_{m},\sigma_{m}) \mid \mu_{m} \in \real ,\, \sigma_{m} \in (0,\infty) \}$ denote the parameter space of the $m$-th component. 
Then the entire parameter space $\Theta$ can be represented as 
\begin{equation}
 \Theta \equiv \{(\alpha_{1},\dots,\alpha_{m}) \mid \sum_{m=1}^{M}\alpha_{m}=1 ,\, \alpha_{m}\geq 0\} \times \prod_{m=1}^{M} \Omega_{m}. 
 \nonumber 
\end{equation}

For a given sample $\vr{X} \equiv (X_{1},\dots,X_{n})$ from $f(x;\theta_{0})$, the likelihood function is defined as 
\begin{equation}
 l(\theta;\vr{X}) \equiv \prod_{i=1}^{n} f(X_{i};\theta) = \prod_{i=1}^{n} \left\{\sum_{m=1}^{M}\alpha_{m}f_{m}(X_{i};\mu_{m},\sigma_{m})\right\}. 
  \nonumber 
\end{equation}

Throughout this paper, we fix $M$, the number of components of mixture models. 
Let $\scrf_{m} \subset \{f(x;\theta) \mid \theta \in \Theta \}$ denote the set of location-scale mixture densities which consist of no more than $m$ components. 
For example, if a mixture density satisfies $\alpha_{m+1}=\dots=\alpha_{M}=0$, then the density belongs to $\scrf_{m}$. 
Note that $\scrf_{M} = \{f(x;\theta) \mid \theta \in \Theta\}$. 

Let $\sigma_{(1)}$ and $\sigma_{(M)}$ denote the minimum and the maximum values of the scale parameters: 
\begin{equation}
 \sigma_{(1)} \equiv \min_{1 \leq m \leq M} \sigma_{m} 
  \quad ,\quad 
  \sigma_{(M)} \equiv \max_{1 \leq m \leq M} \sigma_{m} \; . 
  \label{eq:orderedSigma} 
\end{equation}
Let $\{c_{n}\}_{n=1}^{\infty}$ and $\{b_{n}\}_{n=1}^{\infty}$ denote sequences of positive reals which converge to zero. 
In our discussion, we use two constrained parameter space. 
Define $\Theta_{c_{n}}, \Theta_{b_{n}}$ as follows: 
\begin{eqnarray}
  \Theta_{c_{n}} \equiv \{\theta \in \Theta \mid \sigma_{(1)} \geq c_{n}\}, 
   \quad 
  \Theta_{b_{n}} \equiv \{\theta \in \Theta \mid \frac{\sigma_{(1)}}{\sigma_{(M)}} \geq b_{n}\}. 
  \nonumber 
\end{eqnarray}

\subsection{Regularity conditions}

We introduce assumptions for the strong consistency of the maximum likelihood estimator. 
These assumptions are essentially the same as in \citet{W1949}, \citet{R1981} and \citet{TT2003-35}. 

Let $\Gamma$ denote any compact subset of $\Theta$. 
\begin{assumption}
 For $\theta \in \Theta$ and any positive real number $\rho$, let 
 \begin{eqnarray}
  f(x;\theta,\rho) 
   & \equiv & \sup_{\dist(\theta',\theta) \leq \rho}f(x;\theta')
   , 
   \nonumber
 \end{eqnarray}
 where $\dist(\theta',\theta)$ is a distance between $\theta$ and $\theta'$. 
 For each $\theta \in \Gamma$ and sufficiently small $\rho$, 
 $f(x;\theta, \rho)$ is measurable. 
\label{assumption:1}
\end{assumption}
\begin{assumption}
 For each $\theta \in \Gamma$,  if\/
 $\lim_{j \rightarrow \infty}\theta^{(j)} = \theta, \; (\theta^{(j)} \in \Gamma)$ 
 then\/
 $\lim_{j \rightarrow \infty}
 f(x;\theta^{(j)}) = f(x;\theta)$
 except on a set which is a null set and 
 does not depend on the sequence 
 $\{\theta^{(j)}\}_{j=1}^{\infty}$.
\label{assumption:2}
\end{assumption}
\begin{assumption}
 \begin{eqnarray}
  \int \abslr{
   \log{f(x;\theta_0)}
   }
   f(x;\theta_0)\rmd x < \infty .
   \nonumber 
 \end{eqnarray}
\label{assumption:3}
\end{assumption}

Furthermore, in Section~\ref{sec:main-results}, we impose Assumption~\ref{assumption:4} or \ref{assumption:5} 
according to what type of penalty is made. 
If the penalty is imposed on the scale parameters themselves, then we impose Assumption~\ref{assumption:4}. 
Alternatively, if the penalty is imposed on the ratios of the scale parameters, then we impose Assumption~\ref{assumption:5}. 
\begin{assumption}
 There exist real constants 
 $v_{0},v_{1} > 0$ and $\beta > 1$ such that 
 \begin{eqnarray}
  f_{m}(x;\mu_{m}=0, \sigma_{m}=1) 
   \leq 
   \min \{v_{0} \;, \; v_{1} \cdot \abs{x}^{-\beta}\}
   \nonumber 
 \end{eqnarray}
 for all $m$.
\label{assumption:4}
\end{assumption}
\begin{assumption}
 There exist real constants 
 $v_{0},v_{1} > 0$ and $\beta > 2$ such that 
 \begin{eqnarray}
  f_{m}(x;\mu_{m}=0, \sigma_{m}=1) 
   \leq 
   \min \{v_{0} \;, \; v_{1} \cdot \abs{x}^{-\beta}\}
   \nonumber 
 \end{eqnarray}
 for all $m$.
\label{assumption:5}
\end{assumption}
Note that Assumption~\ref{assumption:5} is stronger than Assumption~\ref{assumption:4}. 
Therefore, if Assumption~\ref{assumption:1},\ref{assumption:2},\ref{assumption:3} and \ref{assumption:5} hold, 
then Assumption~\ref{assumption:1},\ref{assumption:2},\ref{assumption:3} and \ref{assumption:4} hold.

\subsection{Strong Consistency}

According to \citet{R1981}, we define strong consistency of estimators for mixture distributions 
by identifying the parameters whose densities are equal. 
Let 
\begin{equation}
 \Theta(\theta') \equiv \{\theta \in \Theta \mid f(x;\theta) = f(x;\theta') \text{ almost everywhere on } \real \} 
  \nonumber 
\end{equation}
Furthermore we abbreviate $\Theta(\theta_{0})$ as $\Theta_{0}$. 
Given $U , V \subset \Theta$, the distance between $U$ and $V$ are defined as 
\begin{equation}
 \dist(U, V) \equiv \inf_{\theta \in U}\inf_{\theta' \in V} \dist(\theta, \theta'). 
 \nonumber 
\end{equation}
We now define strong consistency of estimators for mixture distributions as follows. 
\begin{definition}
 An estimator $\hat{\theta}_{n}$ is strongly consistent iff 
 \begin{equation}
  \prob\left(\lim_{n \rightarrow \infty } \dist(\Theta(\hat{\theta}_{{n}}), \Theta_{0}) = 0 \right) = 1.   
  \nonumber 
 \end{equation}
\end{definition}
In this paper, two notations ``$\prob(A)=1$'' and ``$A, a.s.$''
($A$ holds almost surely), will be used interchangeably.

\section{Main results}
\label{sec:main-results}

\subsection{Consistency of penalized maximum likelihood estimator 
when the penalty is imposed on the minimum of the ratios of the scale parameters}
\label{sec:penal-ratio-scale}

Now we define a penalized likelihood. 
Let $\bar{r}_{n}(\cdot)$ denote a function on $(0, 1]$ which satisfies the following assumption 
and is not identically equal to zero. 
\begin{assumption} 
 \label{assumption:6}
 $\exists \bar{R} < \infty , \; \exists \bar{r} > 0 , \; \exists \delta > 0 , \; 0 \leq \exists \bnd < 1$ such that 
 \begin{equation}
    \quad 0 \leq \bar{r}_{n}(y) \leq \min{\{ \bar{R} , \; \bar{r} \cdot y^{M+\delta} \cdot \exp{(n^{\bnd})}\}}. 
   \nonumber
 \end{equation} 
\end{assumption}
The Assumption~\ref{assumption:6} means that 
$\bar{r}_{n}(y)$ is nonnegateve, bounded in $n$ and $y$, and converges to zero sufficiently fast as $y$ tends to zero. 
Note that we can take a discontinuous function as $\bar{r}_{n}(y)$. 
In Corollary~\ref{cor:const-ratio-scale}, 
we obtain the consistency of a constrained maximum likelihood estimator 
by using a discontinuous penalty function. 

We define a penalty function $1/r_{n}(\theta)$ or a reward function $r_{n}(\theta)$ as 
\begin{equation}
 r_{n}(\theta) \equiv \bar{r}_{n}\left(\frac{\sigma_{(1)}}{\sigma_{(M)}}\right). 
  \nonumber
\end{equation}
The penalized likelihood function is defined as 
$g_{n}(\theta;\vr{X}) \equiv l(\theta;\vr{X}) \cdot r_{n}(\theta)$. 
The penalized maximum likelihood estimator is defined as 
$\hat{\theta}_{g_{n}} \equiv \argsup_{\theta \in \Theta} g_{n}(\theta;\vr{X})$. 
As stated in Section~\ref{sec:preliminaries}, 
the likelihood $l(\theta;\vr{X})$ may increase to infinity as $\sigma_{m}$ decreases to zero. 
However, if the penalty $1/r_{n}(\theta)$ increases to infinity or $r_{n}(\theta)$ decreases to zero, 
the divergence of the likelihood may be avoided. 
This happens when a part of the scale parameters decreases to zero. 
If all the scale parameters decreases to zero, 
then the likelihood $l(\theta;\vr{X})$ decreases to zero 
because a component with a very small scale parameter can only have a large contribution to a single observation at most. 
Therefore, the existence of the penalty term may prevent the positive divergence of the likelihood. 

Let $b_{0} > 0$. 
In this section, 
we take $b_{n}$ as follows: 
\begin{equation}
  b_{n} = b_{0}\cdot \exp{(-n^{\bnd})}
 \nonumber 
\end{equation}
We also assume the following conditions. 
\begin{assumption}
 \label{assumption:8}
 There exist a positive real $r_{\Theta_{0}}$ and a positive integer $N$ such that 
 $r_{n}(\theta) \geq r_{\Theta_{0}}$
 for all $\theta \in \Theta_{0}$ and $n \geq N$. 
\end{assumption}
If $\bar{r}_{n}(y)$ is positive and unimodal or increasing, then $r_{n}(\theta)$ satisfies Assumption~\ref{assumption:8}. 
Assumption~\ref{assumption:8} guarantees that the penalized likelihood is nearly unaffected by the penalty term 
for $\theta \in \Theta_{0}$ when sample size $n$ is large. 
\begin{assumption}
 There exist real constants $\cnd$, $c_{0}$ and $\Delta$ such that $0 \leq \bnd < \cnd < 1$, $c_{0} > 0$, $\Delta > 0$
 and the following relation holds for all $\theta \in \Theta_{c_{n}}$ 
 and $n \in \naturalnum$ {\rm :} 
 \begin{equation}
  r_{n}(\theta) > (\sigma_{(1)})^{M}
   \quad \Rightarrow  \quad
   \sigma_{(M)} < \frac{(\sigma_{(1)})^{\Delta}}{b_{n}}, 
   \nonumber 
 \end{equation}
 where $c_{n} = c_{0}\cdot\exp{(-n^{\cnd})}$ and $\Theta_{c_{n}} = \{\theta \in \Theta \mid \sigma_{(1)} \geq c_{n}\}$. 
 \label{assumption:9}
\end{assumption}
Assumption~\ref{assumption:9} means that all the scale parameters of $\theta \in \Theta_{c_{n}}$ are equally small 
if $r_{n}(\theta) > (\sigma_{(1)})^{M}$. 

The assumptions for the penalties are not so restrictive. 
For example, if we set 
$\bar{r}_{n}(y) = \bar{r}\cdot y^{\alpha-1}\cdot e^{n^{\tilde{d}}}$ 
and assume $\alpha > M+1$, 
then $\bar{r}_{n}(y)$ satisfies the Assumption~\ref{assumption:6} and 
$r_{n}(\theta) = \bar{r}_{n}(\frac{\sigma_{(1)}}{\sigma_{(M)}})$ satisfies the Assumption~\ref{assumption:8} and \ref{assumption:9}. 

Then the following theorem holds. 
\begin{thm}
 \label{thm:penalt-ratio-scale}
 Suppose that $\scrf_{M}$ satisfies the Assumption~\ref{assumption:1},\ref{assumption:2},\ref{assumption:3} and \ref{assumption:5}, 
 and $f(x;\theta_{0}) \in \scrf_{M}\setminus\scrf_{M-1}$. 
 Suppose that the penalty function $r_{n}(\theta)$ satisfies the Assumption~\ref{assumption:6}, 
 \ref{assumption:8} and \ref{assumption:9}. 
 Then the penalized maximum likelihood estimator $\hat{\theta}_{g_{n}}$ is strongly consistent. 
\end{thm}

A proof of Theorem~\ref{thm:penalt-ratio-scale} is given in section~\ref{sec:proof-thm:penalt-ratio-scale}. 

As a corollary of Theorem~\ref{thm:penalt-ratio-scale}, 
we can obtain the consistency of a constrained maximum likelihood estimator. 
Let us define the constrained maximum likelihood estimator restricted to $\Theta_{b_{n}}$ as 
\begin{equation}
 \hat{\theta}_{b_{n}} \equiv \argsup_{\theta \in \Theta_{b_{n}}}l(\theta;\vr{X}). 
  \nonumber 
\end{equation}
If we put $\bar{r}_{n}(y)$ and $r_{n}(\theta)$ as 
\begin{equation}
 \bar{r}_{n}(y) = 
  \begin{cases}
   1 & (y \geq b_{n}) \\ 
   0 & (y < b_{n})
  \end{cases}
  \quad , \quad 
  r_{n}(\theta) = \bar{r}_{n}\left(\frac{\sigma_{(1)}}{\sigma_{(M)}}\right) = 
  \begin{cases}
   1 & ({\sigma_{(1)}}/{\sigma_{(M)}} \geq b_{n}) \\ 
   0 & ({\sigma_{(1)}}/{\sigma_{(M)}} < b_{n})
  \end{cases}
  , 
  \label{penaltyAsConstraint:eq}
\end{equation}
then $\hat{\theta}_{b_{n}}$ is equal to the penalized maximum likelihood estimator 
$\hat{\theta}_{g_{n}} = \argsup_{\theta \in \Theta} g_{n}(\theta;\vr{X}) = \argsup_{\theta \in \Theta} l(\theta;\vr{X}) \cdot r_{n}(\theta)$. 
If we take $0 < \bnd < 1$, 
then $r_{n}(\theta)$ given in (\ref{penaltyAsConstraint:eq}) satisfies Assumption~\ref{assumption:6}, 
\ref{assumption:8} and \ref{assumption:9}. 
From this and Theorem~\ref{thm:penalt-ratio-scale}, we obtain the following corollary. 
\begin{cor}
 \label{cor:const-ratio-scale}
 Suppose that $\scrf_{M}$ satisfies the Assumption~\ref{assumption:1},\ref{assumption:2},\ref{assumption:3} and \ref{assumption:5}, 
 and $f(x;\theta_{0}) \in \scrf_{M}\setminus\scrf_{M-1}$. 
 If we take $0 < \bnd < 1$, 
  then the constrained maximum likelihood estimator $\hat{\theta}_{b_{n}}$ 
 is strongly consistent.
\end{cor}
By Corollary \ref{cor:const-ratio-scale}, the problem stated in \citet{H1985} is solved positively.

\subsection{Consistency of penalized maximum likelihood estimator when the penalties are imposed on the scale parameters themselves}
\label{sec:penal-scale-param}

We also give a consistency result for the penalized maximum likelihood estimator in which the penalties are imposed on the scale parameters themselves.  
Let $\bar{s}_{n}(\cdot)$ denote a function on $(0, \infty)$ which satisfies the following assumptions. 
\begin{assumption}

 $\bar{s}_{n}(y)$ is nonnegative, uniformly bounded and not identically equal to zero: 
 \begin{equation}
  0 \leq \bar{s}_{n}(y) \leq \bar{S} < \infty 
   \quad , \quad 
   \sup_{y > 0}\bar{s}_{n}(y) > 0 . 
   \nonumber
 \end{equation}
 \label{assumption:10}
\end{assumption}
\begin{assumption}
 $\bar{s}_{n}(y)$ converges to zero sufficiently fast as $y$ tends to zero: 
 \begin{equation}
   \exists \bar{s} > 0 , \; 0 \leq \exists \cnd < 1 
   \quad \text{ s.t. } \quad 
   0 < \sup_{y > 0} \frac{\bar{s}_{n}(y)}{y^{M}} \leq \bar{s} \cdot \exp{(n^{\cnd})}
   \nonumber
 \end{equation} 
 \label{assumption:11}
\end{assumption}

Then we define a penalty function $1/s_{n}(\theta)$ or reward function $s_{n}(\theta)$ as follows: 
\begin{equation}
 s_{n}(\theta) \equiv \prod_{m=1}^{M} \bar{s}_{n}(\sigma_{m}). 
  \nonumber
\end{equation}
The penalized likelihood function is defined as 
$h_{n}(\theta;\vr{X}) \equiv l(\theta;\vr{X}) \cdot s_{n}(\theta)$. 
The penalized maximum likelihood estimator is defined as 
$\hat{\theta}_{h_{n}} \equiv \argsup_{\theta \in \Theta} h_{n}(\theta;\vr{X})$. 

We also assume the following condition. 
\begin{assumption}
 \label{assumption:12}
 There exist a positive real $s_{\Theta_{0}}$ and a positive integer $N$ such that 
 $s_{n}(\theta) \geq s_{\Theta_{0}}$ 
 for all $\theta \in \Theta_{0}$ and $n \geq N$. 
\end{assumption}

The assumptions for the penalty are not so restrictive. 
For example, if we set $\bar{s}_{n}(y) = {e^{-\frac{\beta}{y}}}\cdot{y^{-(\alpha+1)}}$ and assume $\alpha, \beta > 0$, 
then $\bar{s}_{n}(y)$ satisfies the Assumption~\ref{assumption:10} and \ref{assumption:11}, and 
$s_{n}(\theta) = \prod_{m=1}^{M}\bar{s}_{n}(\sigma_{m})$ satisfies the Assumption~\ref{assumption:12}.

We now state the consistency of the penalized maximum likelihood estimator 
when the penalty is imposed on the scale parameters themselves. 
\begin{thm}
 Suppose that $\scrf_{M}$ satisfies the Assumption~\ref{assumption:1},\ref{assumption:2},\ref{assumption:3} and \ref{assumption:4}, 
 and $f(x;\theta_{0}) \in \scrf_{M}\setminus\scrf_{M-1}$. 
 Suppose that the penalty function $s_{n}(\theta)$ satisfies the Assumption~\ref{assumption:10},\ref{assumption:11} and \ref{assumption:12}. 
 Then the penalized maximum likelihood estimator $\hat{\theta}_{h_{n}}$ is strongly consistent. 
 \label{thm:penalt-scale-param}
\end{thm}

The statement of Theorem \ref{thm:penalt-scale-param} is an extension of Corollary 1 of \citet{CRI2003}. 
In their statement, penalties for the location parameters $\mu_{1},\cdots, \mu_{M}$ may be required. 
This is because, in their proof, they use a compactification of the parameter space, 
but their penalized likelihood is not continuous over the compactified parameter space. 
For example, if $\mu_{1} \rightarrow \infty$, 
then other components may still exist and hence their penalized likelihood may not tend to zero. 

We give a proof of Theorem \ref{thm:penalt-scale-param} in section~\ref{sec:proof-thm:penalt-scale-param}. 

\section{Proofs}
\label{sec:proofs}

In this section, we prove Theorem~\ref{thm:penalt-ratio-scale} and \ref{thm:penalt-scale-param}. 
The organization of this section is as follows. 
In section \ref{sec:some-lemmas}, we state some lemmas needed for proving Theorem~\ref{thm:penalt-ratio-scale} and \ref{thm:penalt-scale-param}. 
Section~\ref{sec:proof-thm:penalt-ratio-scale} and \ref{sec:proof-thm:penalt-scale-param} are devoted to the proof of 
Theorem~\ref{thm:penalt-ratio-scale} and \ref{thm:penalt-scale-param} respectively. 

\subsection{Some lemmas}
\label{sec:some-lemmas}

We state some lemmas needed for proving Theorem~\ref{thm:penalt-ratio-scale} and \ref{thm:penalt-scale-param}. 
Proofs of Lemma~\ref{thm:constr-scale-param}, \ref{lem:importantInterval}, \ref{thm:kappa-lambda}, 
\ref{lem:boundingComponentsWithKappa}, \ref{lem:expectation-equality} and \ref{lem:WaldThm1Type} 
are given in the longer version of \citet{TT2003-35}. 

In \citet{TT2003-35}, we showed that 
when the constraint is appropriately imposed on the minimum of the scale parameters, 
the constrained maximum likelihood estimator is strongly consistent under regularity conditions. 
Let us define the constrained maximum likelihood estimator restricted to $\Theta_{c_{n}} = \{\theta \in \Theta \mid \sigma_{(1)} \geq c_{n}\}$ 
by $\hat{\theta}_{c_{n}} \equiv \argsup_{\theta \in \Theta_{c_{n}}}l(\theta;\vr{X})$. 
\begin{lem}{\rm (\citet{TT2003-35})} 
 \label{thm:constr-scale-param}
 Suppose that $\scrf_{M}$ satisfies the Assumption~\ref{assumption:1},\ref{assumption:2},\ref{assumption:3} and \ref{assumption:4}, 
 and $f(x;\theta_{0}) \in \scrf_{M}\setminus\scrf_{M-1}$. 
 Let $c_{0} > 0$ and $0 < \cnd < 1$. 
 If $c_{n} = c_{0}\cdot \exp(-n^{\cnd})$, 
 then the constrained maximum likelihood estimator 
 $\hat{\theta}_{c_{n}}$ restricted to 
 $\Theta_{c_{n}}$ 
 is strongly consistent.
\end{lem}
As in the case of uniform mixture in \cite{TT2003-20}, it is readily verified that 
if $b_n$ decreases to zero faster than $e^{-n}$, 
then the consistency of the constrained maximum likelihood estimator fails.
Therefore, the rate obtained in Lemma~\ref{thm:constr-scale-param} is almost the lower bound of $b_n$ which maintains the strong consistency. 

Let 
\begin{eqnarray}
 X_{n,1}  \equiv  \min{\{X_{1}, \ldots, X_{n}\}}
  \quad , \quad 
 X_{n,n}  \equiv  \max{\{X_{1}, \ldots, X_{n}\}}.
  \nonumber 
\end{eqnarray}
\begin{lem}
 \label{lem:importantInterval}
 {\rm(\citet{TT2003-35})}
 Suppose that Assumption~\ref{assumption:4} 
 is satisfied. 
 For any real positive constants $A_{0} > 0 , \zeta > 0$, let 
 \begin{equation}
  A_{n} \equiv A_{0}\cdot n^{\frac{2+\zeta}{\beta-1}}, 
   \label{eq:An} 
 \end{equation}
 where $\beta$ is defined by Assumption~\ref{assumption:4}. 
 Then
 \begin{equation}
  \prob\left(
	X_{n,1} < -A_{n} \; \rmor \;  X_{n,n} > A_{n}
	\quad \io
       \right) = 0. 
  \nonumber 
 \end{equation}
\end{lem}
where $\io$ means ``infinitely often''. 
By Lemma~\ref{lem:importantInterval}, we can bound
the behavior of the minimum and the maximum of the sample with probability 1. 
In the following sections, we take $A_{0}$ large 
enough to satisfy (\ref{eq:A-zero:condition:thm:penalt-ratio-scale}) 
and ignore the event $\{X_{n,1} < -A_{n} \; \rmor \;  X_{n,n} > A_{n}\}$. 

Let $R_{n}(V)$ denote the number of observation which belong to a set $V \subset \real$: 
\begin{equation}
 R_{n}(V) \equiv \sharp{\{X_{i} \mid X_{i} \in V,\; i=1,\dots,n\}}. 
 \nonumber 
\end{equation}
Let $P_{0}(V)$ denote the probability of
$V \subset \real$ under the true density: 
\begin{eqnarray}
P_{0}(V) \equiv \int_{V} f(x;\theta_{0}) \rmd x . 
 \nonumber 
\end{eqnarray}

\begin{figure}[htbp]
 \begin{center}
  \includegraphics[width=10cm]{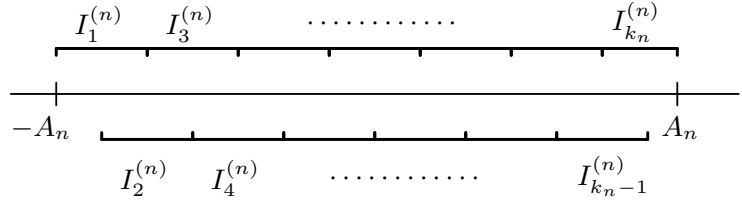}
  \caption{$I_{1}^{(n)},I_{2}^{(n)},\dots,I_{k_{n}}^{(n)}$}
  \label{fig:coverage}
 \end{center}
\end{figure}

Let us consider an interval $[\mu - w_{n},\, \mu + w_{n}]$ with the center $\mu$ and the length $2w_{n}$. 
If $w_{n}=0$, then $R_{n}([\mu - w_{n},\, \mu + w_{n}])$ is clearly 0. 
In the following lemma, we state that 
if $w_{n}$ decreases to zero faster than a power of $1/n$, 
then $R_{n}([\mu - w_{n},\, \mu + w_{n}]) < 2$ holds for every $\mu \in \real$ with probability 1. 
\begin{lem}
 \label{lem:nondense}
 Suppose that Assumption~\ref{assumption:4} 
 is satisfied. 
 Let $\{w_{n}\}_{n=1}^{\infty}$ be a sequence of real numbers which satisfies
 \begin{equation}
  \lim_{n\rightarrow \infty}n^{3+\delta'}\cdot A_{n}\cdot w_{n} = 0, 
   \label{eq:assumption:lem:nondense}
 \end{equation}
 where $\delta' > 0$ and $A_{n}$ is defined by (\ref{eq:An}). 
 Then 
 \begin{equation}
  \prob\left(
	\sup_{\mu \in \real}
	R_{n}( [\mu - {w_{n}},\, \mu + {w_{n}}] ) > 1 
	\quad 
	\io 
       \right) 
  = 0. 
  \nonumber 
 \end{equation}
\end{lem}
\Proof 
From Lemma~\ref{lem:importantInterval}, we ignore the event $\{X_{n,1} < -A_{n} \; \rmor \;  X_{n,n} > A_{n}\}$. 
Then 
\begin{eqnarray}
 \sup_{\mu \in \real}R_{n}( [\mu - {w_{n}},\, \mu + {w_{n}}] )  > 1 
 \quad \Leftrightarrow \quad 
 \sup_{\mu \in [-A_{n}+{w_{n}},\, A_{n}-{w_{n}}]}R_{n}([\mu - {w_{n}},\, \mu + {w_{n}}] )  > 1 
 \quad a.s. \nonumber \\
 \label{eq:equivalenceOfevents:lem:nondense}
\end{eqnarray}
Now we cover $[-A_{n},A_{n}]$ by short intervals of length $4w_{n}$. 
Let 
\begin{eqnarray}
 \lefteqn{
  I_{1}^{(n)} \equiv [-A_{n},\; -A_{n} + 4w_{n}]
  , \; 
  I_{2}^{(n)} \equiv [-A_{n} + 2w_{n},\; -A_{n} + 6w_{n}]\,,\dots\,,
  } & & \nonumber \\ 
 &  & 
  I_{k_{n}-1}^{(n)} \equiv [-A_{n} + (k_{n}-6)\cdot w_{n},\; -A_{n} + (k_{n} - 2)\cdot w_{n} ], 
  \nonumber \\
 & & 
  I_{k_{n}}^{(n)} \equiv [-A_{n} + (k_{n}-4)\cdot w_{n},\; A_{n}], 
  \nonumber 
\end{eqnarray}
where $k_{n} \equiv \min \{ k \in \naturalnum \mid k \cdot (2w_{n})  > 2A_{n} \}$. 
See Figure~\ref{fig:coverage}. 
Since any half-open interval of length $2w_{n}$ in $[-A_{n}, A_{n}]$ is covered by one of $I_{1}^{(n)},\dots,I_{k_{n}}^{(n)}$, 
the following relation holds. 
\begin{eqnarray}
 \sup_{\mu \in [-A_{n}+{w_{n}},\, A_{n}-{w_{n}}]}R_{n}([\mu - {w_{n}},\, \mu + {w_{n}}] )  > 1 
  \quad \Rightarrow \quad 
  1 \leq \exists k \leq k_{n} \; , \; R_{n}(I_{k}^{(n)})  > 1 
  \qquad
  \label{eq:inclusionRelationOfevents:lem:nondense}
\end{eqnarray}
Let $u_{0} \equiv \sup_{x}f(x;\theta_{0})$. 
Because $R_n(I_{k}^{(n)}) \sim \bin(n , P_{0}(I_{k}^{(n)}))$ and $P_{0}(I_{k}^{(n)}) \le 2w_{n} u_{0}$, 
we obtain 
\begin{eqnarray}
 \lefteqn{
  \prob
  \left(
  1 \leq \exists k \leq k_{n} \; , \; R_{n}(I_{k}^{(n)})  > 1 
  \right)
  \leq
  \sum_{k=1}^{k_{n}}
  \prob
  \left(
   R_n(I_{k}^{(n)})  > 1 
	 \right)
  } \hspace{2cm} & & 
\nonumber \\ & \leq &
 k_{n}
 \cdot 
 \left\{
 \max_{1\le k \leq k_{n}}
 \prob(R_n(I_{k}^{(n)})  > 1 )
 \right\}
\nonumber \\ & \leq & 
 \left(
  \frac{A_{n}}{w_{n}} + 1
 \right) \cdot 
 \sum_{k=2}^{n}
 \binom{n}{k}
 (2w_{n} u_{0})^{k}(1 - 2w_{n} u_{0})^{n-k}
 \nonumber \\ & \leq &
 \left(
  \frac{A_{n}}{w_{n}} + 1
 \right) \cdot 
 \sum_{k=2}^{n}\frac{n^{k}}{k!}(2w_{n} u_{0})^{k}
 \nonumber \\ & \leq &
   \left(
  \frac{A_{n}}{w_{n}} + 1
 \right) \cdot 
 (2nw_{n} u_{0})^{2} \cdot 
 \exp{(2nw_{n} u_{0})}
 \; .
 \label{eq:upperBoundOfProb:proof:lem:nondense}
\end{eqnarray}
From (\ref{eq:assumption:lem:nondense}), 
when we sum the right hand side of (\ref{eq:upperBoundOfProb:proof:lem:nondense}) over $n$, the resulting series converges. 
Hence by (\ref{eq:equivalenceOfevents:lem:nondense}), (\ref{eq:inclusionRelationOfevents:lem:nondense}), (\ref{eq:upperBoundOfProb:proof:lem:nondense}) 
and Borel-Cantelli lemma, we have 
\begin{equation}
 \prob\left(
       \sup_{\mu \in \real}
       R_{n}([\mu - {w_{n}},\, \mu + {w_{n}}])  > 1  
       \quad 
       \io 
      \right) 
 = 0. 
 \nonumber 
\end{equation}
\qed

\begin{figure}[htbp]
 \begin{center}
  \includegraphics[width=8cm]{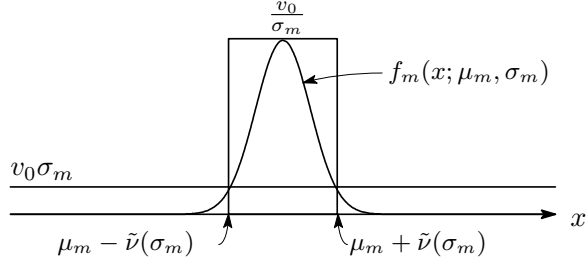}
  \caption{Each component is bounded by a step function.}
  \label{fig:boundedByStepFunction}
 \end{center}
\end{figure}
Next we bound the component densities from above. 
For $\beta>2$, define $\tilde{\nu}(\sigma)$ as 
\begin{equation}
 \tilde{\nu}(\sigma) \equiv \left(\frac{v_{1}}{v_{0}}\right)^{\frac{1}{\beta}}\cdot\sigma^{1-\frac{2}{\beta}}. 
 \label{eq:definitionOfTildeNu} 
\end{equation}
Let $\indicator_{U}(x)$ denote the indicator function of $U \subset \real$.
\begin{lem}
 \label{lem:boundingComponentsWithSmallSigma}
 Suppose that Assumption~\ref{assumption:5} 
 is satisfied. 
 Then the following inequalities hold. 
 \begin{equation}
   f_{m}(x;\mu_{m}, \sigma_{m}) \leq \max\{1_{[\mu_{m}-\tilde{\nu}(\sigma_{m}),\; \mu_{m}+\tilde{\nu}(\sigma_{m})]}(x)\cdot\frac{v_{0}}{\sigma_{(1)}}
   \; , \; {v_{0}}{\sigma_{(M)}}\}
   \; , \quad 
   1 \leq m \leq M
   . 
   \label{eq:boundingComponentsWithSmallSigma:lem:boundingComponentsWithSmallSigma}
 \end{equation}
\end{lem}
\Proof 
From Assumption~\ref{assumption:5}, each component is bounded from above as 
\begin{equation}
 f_{m}(x;\mu_{m}, \sigma_{m}) \leq \max\{1_{[\mu_{m}-\tilde{\nu}(\sigma_{m}),\; \mu_{m}+\tilde{\nu}(\sigma_{m})]}(x)\cdot\frac{v_{0}}{\sigma_{m}}
  \; , \; {v_{0}}{\sigma_{m}}\}. 
  \nonumber 
\end{equation}
See figure~\ref{fig:boundedByStepFunction}. 
From this and (\ref{eq:orderedSigma}), we obtain (\ref{eq:boundingComponentsWithSmallSigma:lem:boundingComponentsWithSmallSigma}). 
\qed 

Let 
$E_{0}[\cdot]$ denote the expectation under the true parameter $\theta_{0}$.
\begin{lem} {\rm(\citet{TT2003-35})}
 \label{thm:kappa-lambda}
 Suppose that $\scrf_{M}$ satisfies the Assumption \ref{assumption:1},\ref{assumption:2},\ref{assumption:3} and \ref{assumption:4},  
 and $f(x;\theta_{0}) \in \scrf_{M} \backslash \scrf_{M-1}$. 
Then there  exist real constants $\kappa, \lambda > 0$ such that 
 \begin{eqnarray}
  E_{0}
  \left[
   \log{
   \left\{ 
    f(x;\theta) + \kappa
   \right\}
   }
  \right]
  +
  \lambda
  <
  E_{0}[\log{f(x;\theta_{0})}] 
  \label{eq:kappa-lambda}
 \end{eqnarray}
 for all 
 $f(x;\theta) \in \scrf_{M-1}$. 
\end{lem} 

Fix arbitrary $\kappa_{0}> 0$, which corresponds to $\kappa$ in Lemma~\ref{thm:kappa-lambda}. 
For $\beta > 1$, define $\nu(\sigma)$ as 
\begin{equation}
 \nu(\sigma) \equiv \left(\frac{v_{1}}{\kappa_{0}}\right)^{\frac{1}{\beta}}\cdot\sigma^{1-\frac{1}{\beta}}. 
  \nonumber 
\end{equation}
In a manner similar to the proof of Lemma~\ref{lem:boundingComponentsWithSmallSigma}, we can show the following lemma. 
\begin{lem} {\rm(\citet{TT2003-35})} 
\label{lem:boundingComponentsWithKappa}
Suppose that Assumption~\ref{assumption:4} is satisfied. 
Then the following inequality holds. 
\begin{eqnarray}
 f_{m}(x; \mu_{m}, \sigma_{m}) 
  \leq 
  \max\{\indicator_{[\mu_{m}-\nu(\sigma_{m}),\; \mu_{m}+\nu(\sigma_{m})]}(x) \cdot \frac{v_0}{\sigma_{m}}
  \; , \; {\kappa_{0}}\}. 
  \nonumber 
\end{eqnarray}
\end{lem}

Lemma~\ref{lem:boundingComponentsWithKappa} bounds the tails of a density in a different way than in Lemma~\ref{lem:boundingComponentsWithSmallSigma}. 
On the one hand,  in Lemma~\ref{lem:boundingComponentsWithSmallSigma}, the tails of a density is bounded by the value of scale parameter 
and Assumption~\ref{assumption:5} is needed because $\beta$ should be larger than 2. 
On the other hand, in Lemma~\ref{lem:boundingComponentsWithKappa}, the tails of a density is bounded by a constant 
and only Assumption~\ref{assumption:4} is needed. 
Lemma~\ref{lem:boundingComponentsWithSmallSigma} will be used to prove Theorem~\ref{thm:penalt-ratio-scale}. 
Lemma~\ref{lem:boundingComponentsWithKappa} will be used to prove Theorem~\ref{thm:penalt-ratio-scale} and Theorem~\ref{thm:penalt-scale-param}. 
Therefore, Theorem~\ref{thm:penalt-ratio-scale} needs Assumption~\ref{assumption:5} which is stronger than Assumption~\ref{assumption:4}. 

Let $\scrm$ be a subset of $\{1,2,\ldots,M\}$ and 
let $\abslr{\scrm}$ denote the number of elements in
$\scrm$. 
Denote by $\theta_{\scrm}$ a subvector of $\theta \in \Theta$ 
consisting of the components in $\scrm$. 
Then the parameter space of subprobability measures 
consisting of the components in $\scrm$ 
is
\begin{equation}
 \bar{\Theta}_{\scrm} \equiv 
  \{
  \theta_{\scrm} \mid \theta \in \Theta, \sum_{m \in \scrm} \alpha_m
  \le 1 
  \}.
  \nonumber 
\end{equation}
Corresponding density and the set of subprobability densities are
denoted by 
\begin{align}
 f_{\scrm}(x;\theta_{\scrm})
  &\equiv \sum_{k \in \scrm} \alpha_{k} f_k(x;\mu_k, \sigma_{k}) ,
  \nonumber 
 \\
 \scrf_{\scrm} 
  & \equiv 
  \{
  f_{\scrm}(x;\theta_{\scrm}) 
  \mid 
  \theta_{\scrm} \in \bar{\Theta}_{\scrm}
  \} .
  \nonumber 
\end{align}
Then $\scrf_{K}$, the set of subprobability densities with no more than $K$ components, 
can be represented as  
\begin{equation}
 \scrf_{K} \equiv \bigcup_{\abslr{\scrm}\leq K}\scrf_{\scrm}
  \qquad (1 \leq K \leq M). 
  \nonumber 
\end{equation}

The following lemma follows from the bounded convergence theorem.
\begin{lem} {\rm(\citet{TT2003-35})}
 \label{lem:expectation-equality}
 Let 
 ${\Gamma}_{\scrm}$ 
 denote any compact subset of $ \bar{\Theta}_{\scrm}$.
 For any real constant $\kappa_{0} \geq 0$ and 
 any point $\theta_{\scrm} \in {\Gamma}_{\scrm}$, 
 the following equality holds 
 under Assumption~\ref{assumption:1} and \ref{assumption:3}.
 \begin{eqnarray}
  \lim_{\rho \rightarrow 0}
   E_{0}[\log\{f_{\scrm}(x;\theta_{\scrm}, \rho) + \kappa_{0}\}]
    = 
    E_{0}[\log\{f_{\scrm}(x;\theta_{\scrm}) + \kappa_{0}\}] \ .
    \nonumber 
 \end{eqnarray}
\end{lem}

The following lemma follows from lemma \ref{lem:expectation-equality}. 
\begin{lem} {\rm(\citet{TT2003-35})}
\label{lem:WaldThm1Type}
 Let $\kappa_{0}$ and $\lambda_{0}$ be real constants 
 which corresponds to $\kappa$ and $\lambda$ in Lemma~\ref{thm:kappa-lambda}. 
 Let 
 ${\Gamma}_{\scrm}$ 
 denote any compact subset of $\bar{\Theta}_{\scrm}$.
 Let $\scrb(\theta_{\scrm},\rho(\theta_{\scrm}) )$ denote the open ball
 with center $\theta_{\scrm}$ and radius $\rho(\theta_{\scrm})$. 
 Suppose that Assumption~\ref{assumption:1} and \ref{assumption:3} hold. 
 Then ${\Gamma}_{\scrm}$ can be covered by 
 a finite number of balls \\
 $\scrb(\theta_{\scrm}^{(1)}, \rho(\theta_{\scrm}^{(1)}))
 , \ldots, \scrb(\theta_{\scrm}^{(S)}, \rho(\theta_{\scrm}^{(S)}))
 $
 such that 
\begin{eqnarray}
 E_{0}[\log{\{f_{\scrm}(x;\theta_{\scrm}^{(s)},\rho(\theta_{\scrm}^{(s)})) + \kappa_{0}\}}]
  + \lambda_{0}
  <
  E_{0}[\log{f(x;\theta_0)}]
  \; , \quad (s=1,\ldots,S) \; . 
  \nonumber 
\end{eqnarray}
\end{lem}

\subsection{Proof of Theorem~\ref{thm:penalt-ratio-scale}}
\label{sec:proof-thm:penalt-ratio-scale}

First, we partition the parameter space $\Theta$ into two sets. 
Then the proof of the strong consistency of the penalized maximum likelihood estimator is also partitioned into two parts. 
The proof for one set is obtained immediately by applying the result of Lemma~\ref{thm:constr-scale-param}. 

\subsubsection{Partitioning the parameter space}

Let $\bnd$ be a constant defined by Assumption~\ref{assumption:6}. 
Let $\cnd$ be a constant defined by Assumption~\ref{assumption:9}. 
Define $c_{n} = c_{0}\cdot\exp(-n^{\cnd})$ and $\Theta_{c_{n}}=\{\theta \in \Theta \mid {\sigma_{(1)}} \geq c_{n}\}$. 
Then the parameter space $\Theta$ is divided into two sets: 
\begin{equation}
 \Theta = \Theta_{c_{n}} \cup \Theta_{c_{n}}^{C}, 
  \nonumber 
\end{equation}
where $\Theta_{c_{n}}^{C} = \{\theta \in \Theta \mid \sigma_{(1)} < c_{n}\} $ is the complement of $\Theta_{c_{n}}$. 
From Assumption~\ref{assumption:6}, the reward term $r_{n}(\theta)$ is bounded. 
Furthermore, Assumption~\ref{assumption:8} indicates that 
the asymptotic behavior is not affected by the penalty term around $\Theta_{0}$. 
Therefore the penalized maximum likelihood estimator over $\Theta_{c_{n}}$ is strongly consistent by Lemma~\ref{thm:constr-scale-param}. 
If the maximum of the likelihood function over $\Theta_{c_{n}}^{C}$ is very small, 
then the penalized maximum likelihood estimator over the whole parameter space $\Theta$ is strongly consistent. 
This takes care of $\Theta_{c_{n}}$ and from now on we consider the behavior of the penalized likelihood over $\Theta_{c_{n}}^{C}$. 

Furthermore, we divide $\Theta_{c_{n}}^{C}$ into two sets: 
\begin{eqnarray}
 \Theta_{c_{n}}^{C} & = &  \Phi_{n} \cup \Psi_{n}, 
  \nonumber 
\end{eqnarray}
where
\begin{eqnarray}
 \Phi_{n} & \equiv & \{\theta \in \Theta_{c_{n}}^{C} \mid \frac{r_{n}(\theta)}{(\sigma_{(1)})^{M}} > 1\},  
  \label{eq:definitionOfPhi:thm:penalt-ratio-scale}
 \\ 
 \Psi_{n} & \equiv & \{\theta \in \Theta_{c_{n}}^{C} \mid \frac{r_{n}(\theta)}{(\sigma_{(1)})^{M}} \leq 1\}. 
  \label{eq:definitionOfPsi:thm:penalt-ratio-scale}
\end{eqnarray}
For $\theta \in \Phi_{n}$, all the scale parameters are very small. 
On the other hand, $\theta \in \Psi_{n}$, the penalty $1/r_{n}(\theta)$ is very large and has large contribution relative to the likelihood. 
Therefore, intuitively, it seems that the maximum of the likelihood function over $\Theta_{c_{n}}^{C} = \Phi_{n} \cup \Psi_{n}$ is very small. 
We are going to prove that this is true. 

By the argument used in \citet{W1949}, in order to prove Theorem~\ref{thm:penalt-ratio-scale}, 
it suffices to prove the following two equations.  
\begin{eqnarray}
 \lim_{n \rightarrow \infty}
 \frac{\sup_{\theta \in \Phi_{n}}
  \left\{\prod_{i=1}^{n} f(X_{i};\theta) \right\} \cdot r_{n}(\theta)}
  {\left\{ \prod_{i=1}^{n} f(X_{i};\theta_{0}) \right\} \cdot r_{n}(\theta_{0}) } = 0, \quad a.s. 
  \label{eq:goal-Phi:thm:penalt-ratio-scale}
  \\ 
 \lim_{n \rightarrow \infty}
 \frac{\sup_{\theta \in \Psi_{n}}
  \left\{\prod_{i=1}^{n} f(X_{i};\theta) \right\} \cdot r_{n}(\theta)}
  {\left\{ \prod_{i=1}^{n} f(X_{i};\theta_{0}) \right\} \cdot r_{n}(\theta_{0}) } = 0, \quad a.s. 
  \label{eq:goal-Psi:thm:penalt-ratio-scale}
\end{eqnarray}

\subsubsection{Proof of~(\ref{eq:goal-Phi:thm:penalt-ratio-scale}) for $\Phi_{n}$}
\label{sec:proof:eq:goal-Phi:thm:penalt-ratio-scale}

By the law of large numbers, we have 
\begin{equation}
 \lim_{n\rightarrow\infty} \frac{1}{n} \sum_{i=1}^{n} 
  \log{f(X_{i};\theta_0)} =
  E_0[\log f(x;\theta_0)], \quad a.s.
 \nonumber 
\end{equation}
Furthermore, by Assumption~\ref{assumption:6} and \ref{assumption:8}, we obtain
\begin{equation}
 \lim_{n\rightarrow\infty}\frac{1}{n}\log r_{n}(\theta_{0}) = 0. 
 \nonumber 
\end{equation}
Therefore (\ref{eq:goal-Phi:thm:penalt-ratio-scale}) is implied by 
\begin{eqnarray}
 \limsup_{n\rightarrow \infty} \frac{1}{n} \cdot \sup_{\theta \in \Phi_{n}}
  \left\{ \sum_{i=1}^{n}\log{f(X_{i};\theta)} + \log{r_{n}(\theta)} \right\}
  < E_{0}[\log{f(x;\theta_{0})}] \quad a.s. 
  \label{eq:goal2-Phi:thm:penalt-ratio-scale} 
\end{eqnarray}
Consequently, in order to prove (\ref{eq:goal-Phi:thm:penalt-ratio-scale}), 
it suffices to prove (\ref{eq:goal2-Phi:thm:penalt-ratio-scale}). 

From Assumption~\ref{assumption:9} and (\ref{eq:definitionOfPsi:thm:penalt-ratio-scale}), we have 
\begin{equation}
 \sigma_{(1)} \leq \sigma_{(M)} < \frac{(\sigma_{(1)})^{\Delta}}{b_{n}}
  \quad ,\quad 
  \theta \in \Phi_{n}, 
  \label{eq:inequalitiesForSigmaM:thm:penalt-ratio-scale}
\end{equation}
where $b_{n} = b_{0}\cdot\exp(-n^{\bnd})$
and the first inequality is derived from (\ref{eq:orderedSigma}). 
Because $\sigma_{(1)} < c_{n} = \exp{(-n^{\cnd})}$, we obtain
\begin{equation}
 \sigma_{m} < \exp{(n^{\bnd}-{\Delta\cdot n^{\cnd}})}
  \quad ,\quad 
  1 \leq m \leq M
  \, , \,
  \theta \in \Phi_{n}. 
  \label{eq:inequalitiesForSigma:thm:penalt-ratio-scale}
\end{equation}
Note that $0 \leq \bnd < \cnd < 1 ,\; \Delta > 0$ by Assumption~\ref{assumption:9}. 
Define 
\begin{equation}
 \tilde{J}(\theta) \equiv \bigcup_{m=1}^{M}[\mu_{m}-\tilde{\nu}(\sigma_{m}), \; \mu_{m}+\tilde{\nu}(\sigma_{m})]. 
  \label{eq:definitionOfJ:thm:penalt-ratio-scale}
\end{equation}
Then the following lemma holds. 
\begin{lem}
 \label{lem:numberOfObservationsWithinShortInterval-1:thm:penalt-ratio-scale}
 \begin{equation}
  \prob\left(
	\sup_{\theta \in \Phi_{n}}
	R_{n}(\tilde{J}(\theta)) > M \quad \io
       \right) = 0. 
 \label{eq:lem:numberOfObservationsWithinShortInterval-1:thm:penalt-ratio-scale}
 \end{equation}
\end{lem}
\Proof 
We prove Lemma~\ref{lem:numberOfObservationsWithinShortInterval-1:thm:penalt-ratio-scale} 
by using Lemma~\ref{lem:nondense}. 
Let $w_{n} = \tilde{\nu}(\exp{(n^{\bnd}-{\Delta\cdot n^{\cnd}})})$. 
Because (\ref{eq:definitionOfTildeNu}) and $\beta > 2$ by Assumption~\ref{assumption:5}, 
the assumption (\ref{eq:assumption:lem:nondense}) of Lemma~\ref{lem:nondense} is satisfied. 
From (\ref{eq:inequalitiesForSigma:thm:penalt-ratio-scale}) and (\ref{eq:definitionOfJ:thm:penalt-ratio-scale}), we have 
\begin{eqnarray} 
  \sup_{\theta \in \Phi_{n}} R_{n}(\tilde{J}(\theta)) > M 
 & \Rightarrow & 
  \sup_{\mu \in \real}R_{n}([\mu - w_{n},\, \mu + w_{n}] )  > 1 
  \nonumber 
\end{eqnarray} 
Therefore, by Lemma~\ref{lem:nondense}, we obtain (\ref{eq:lem:numberOfObservationsWithinShortInterval-1:thm:penalt-ratio-scale}). 
\qed 

We now state the following inequality, 
in order to bound the likelihood. 
\begin{lem} For $\theta \in \Phi_{n}$, 
 \label{lem:boundingLikelihoodFunc:thm:penalt-ratio-scale}
 \begin{eqnarray} 
   \sum_{i=1}^{n}\log{f(X_{i};\theta)} 
   \leq 
   R_{n}(\tilde{J}(\theta)) \cdot \log{\frac{v_{0}}{\sigma_{(1)}}} 
  + 
   R_{n}(\tilde{J}(\theta)^{C}) \cdot \log{{v_{0}}{\sigma_{(M)}}}. 
   \nonumber 
 \end{eqnarray}
\end{lem}
\Proof
From Lemma~\ref{lem:boundingComponentsWithSmallSigma} and (\ref{eq:orderedSigma}), 
for $\theta \in \Phi_{n}$, we obtain 
\begin{eqnarray}
 \lefteqn{
 \sum_{i=1}^{n}\log{f(X_{i};\theta)}
   = 
  \sum_{i=1}^{n}\log{\left\{\sum_{m=1}^{M}\alpha_{m}f_{m}(X_{i};\mu_{m}, \sigma_{m})\right\}}
 } & & 
 \nonumber \\ 
 & \leq &
 \sum_{i=1}^{n}\max_{m=1,\dots,M}\log{f_{m}(X_{i};\mu_{m}, \sigma_{m})}
 \nonumber \\ & \leq &
 \sum_{i=1}^{n}
 \max_{m=1,\dots,M} \left\{\max\{
 1_{[\mu_{m}-\tilde{\nu}(\sigma_{m}), \; \mu_{m}+\tilde{\nu}(\sigma_{m})]}(x) \cdot \log{\frac{v_{0}}{\sigma_{(1)}}}  
 \; , \; 
 \log{v_{0}\sigma_{(M)}}
 \}\right\}
 \nonumber \\ & = & 
 R_{n}(\tilde{J}(\theta))\cdot\log{\frac{v_{0}}{\sigma_{(1)}}} + R_{n}(\tilde{J}(\theta)^{C})\cdot\log{v_{0}\sigma_{(M)}}. 
 \nonumber 
\end{eqnarray}
\qed

By Lemma~\ref{lem:boundingLikelihoodFunc:thm:penalt-ratio-scale} and Assumption~\ref{assumption:6}, we obtain for $\theta \in \Phi_{n}$
\begin{eqnarray}
   \sum_{i=1}^{n}\log{f(X_{i};\theta)} + \log{r_{n}(\theta)} 
  \leq 
 R_{n}(\tilde{J}(\theta))\cdot\log{\frac{v_{0}}{\sigma_{(1)}}} 
 + R_{n}(\tilde{J}(\theta)^{C})\cdot\log{v_{0}\sigma_{(M)}} + \log{\bar{R}}. 
  \nonumber 
\end{eqnarray}
Furthermore, from (\ref{eq:inequalitiesForSigmaM:thm:penalt-ratio-scale}), we have for $\theta \in \Phi_{n}$
\begin{eqnarray}
    \sum_{i=1}^{n}\log{f(X_{i};\theta)} + \log{r_{n}(\theta)} 
  \leq 
 R_{n}(\tilde{J}(\theta))\cdot\log{\frac{v_{0}}{\sigma_{(1)}}} 
 + R_{n}(\tilde{J}(\theta)^{C})\cdot\log{\frac{v_{0}(\sigma_{(1)})^{\Delta}}{b_{n}}} + \log{\bar{R}}
 \nonumber \\
 =
 \left(\Delta \cdot R_{n}(\tilde{J}(\theta)^{C}) - R_{n}(\tilde{J}(\theta))\right)\cdot\log{\sigma_{(1)}} 
 - R_{n}(\tilde{J}(\theta)^{C})\cdot\log{b_{n}} + n\log{v_{0}} + \log{\bar{R}}. 
  \nonumber 
\end{eqnarray}
Because $b_{n} = b_{0}\cdot e^{-n^{\tilde{d}}}$, we obtain for $\theta \in \Phi_{n}$
\begin{eqnarray}
  \lefteqn{
    \sum_{i=1}^{n}\log{f(X_{i};\theta)} + \log{r_{n}(\theta)} 
  } & & \nonumber \\
 & \leq &
 \left(\Delta \cdot R_{n}(\tilde{J}(\theta)^{C}) - R_{n}(\tilde{J}(\theta))\right)\cdot\log{\sigma_{(1)}} 
 + R_{n}(\tilde{J}(\theta)^{C})\cdot (n^{\tilde{d}} - \log{b_{0}})
 + n\log{v_{0}} + \log{\bar{R}}
 \nonumber \\
 & \leq &
 \left(\Delta \cdot R_{n}(\tilde{J}(\theta)^{C}) - R_{n}(\tilde{J}(\theta))\right)\cdot\log{\sigma_{(1)}} 
 + n \cdot (n^{\tilde{d}} + \abs{-\log{b_{0}}} + \log{v_{0}}) + \log{\bar{R}}. 
 \nonumber \\
 \label{eq:boundingPenalizedLikelihoodFunc:thm:penalt-ratio-scale}
\end{eqnarray}
By Lemma~\ref{lem:numberOfObservationsWithinShortInterval-1:thm:penalt-ratio-scale}, we obtain
\begin{eqnarray}
  1 
  & = &
  \prob
  \left(
   \bigcup_{N=1}^{\infty}\bigcap_{n=N}^{\infty} \sup_{\theta \in \Phi_{n}}R_{n}(\tilde{J}(\theta)) \leq M
  \right)
   \nonumber \\
   & = &
    \prob
    \left(
     \bigcup_{N=1}^{\infty}\bigcap_{n=N}^{\infty}
     \left\{
      \left\{
       \sup_{\theta \in \Phi_{n}}R_{n}(\tilde{J}(\theta)) \leq M
      \right\}
      \bigcap
      \left\{
       \sup_{\theta \in \Phi_{n}}R_{n}(\tilde{J}(\theta)^{C}) \geq n-M
      \right\}
     \right\}
    \right)
    \nonumber \\
 & \leq &
    \prob
    \left(
     \bigcup_{N=1}^{\infty}\bigcap_{n=N}^{\infty}\sup_{\theta \in \Phi_{n}}
     \left(\Delta \cdot R_{n}(\tilde{J}(\theta)^{C}) - R_{n}(\tilde{J}(\theta))\right)
     \geq \Delta \cdot (n-M) - M
    \right)
    \leq 1. 
    \nonumber \\
   \label{eq:Rn_inequality:thm:penalt-ratio-scale}
\end{eqnarray}
From (\ref{eq:Rn_inequality:thm:penalt-ratio-scale}), the inequality 
$\sup_{\theta \in \Phi_{n}}\Delta \cdot R_{n}(\tilde{J}(\theta)^{C}) - R_{n}(\tilde{J}(\theta))\geq \Delta \cdot (n-M) - M$ 
holds almost surely except for finite number of $n$. 
Therefore, we ignore the event 
$\sup_{\theta \in \Phi_{n}}\Delta \cdot R_{n}(\tilde{J}(\theta)^{C}) - R_{n}(\tilde{J}(\theta)) < \Delta \cdot (n-M) - M$. 
Because $\sigma_{(1)} \leq c_{n} = c_{0}\cdot e^{-n^{d}}$ and (\ref{eq:boundingPenalizedLikelihoodFunc:thm:penalt-ratio-scale}), 
for all sufficiently large $n$ such that $c_{n} \leq 1$ and $\Delta \cdot (n-M) - M \geq 0$ hold, we have 
\begin{eqnarray}
 \lefteqn{
  \sup_{\theta \in \Phi_{n}}
  \left\{
   \sum_{i=1}^{n}\log{f(X_{i};\theta)} + \log{r_{n}(\theta)} 
  \right\} 
  } & & \nonumber \\
 & \leq &
 (\Delta \cdot (n-M) - M)\cdot(-n^{d} + \log{c_{0}})
 + n \cdot (n^{\tilde{d}} + \abs{-\log{b_{0}}} + \log{v_{0}}) + \log{\bar{R}} 
  \qquad a.s.
 \nonumber 
\end{eqnarray}
From Assumption~\ref{assumption:9}, the first term of the right hand side of the above inequality is the main term 
and diverges to $-\infty$ as $n$ increases. 
Therefore, we obtain (\ref{eq:goal2-Phi:thm:penalt-ratio-scale}): 
\begin{equation}
 \limsup_{n \rightarrow \infty} \frac{1}{n} \cdot \sup_{\theta \in \Phi_{n}} 
  \left\{\sum_{i=1}^{n}\log{f(X_{i};\theta)} + \log{r_{n}(\theta)}\right\} = -\infty 
  \quad a.s. 
  \nonumber 
\end{equation}

\subsubsection{Proof of (\ref{eq:goal-Psi:thm:penalt-ratio-scale}) for $\Psi_{n}$}
\label{sec:goal-Psi:thm:penalt-ratio-scale}

The outline of the proof of (\ref{eq:goal-Psi:thm:penalt-ratio-scale}) is as follows. 
First, we partition the parameter space $\Psi_{n}$ into finite subsets $\Psi_{n, \scrk, s}$ 
depending on the set of some parameter smaller than $c_{n}$. 
Then, by using Lemma \ref{lem:nondense}, we can show that the components with $\sigma_{m} < c_{n}$ 
do not contribute to the likelihood more than $M$ data points 
and the contributions are canceld out by the penalty term. 
Therefore, 
from Lemma \ref{lem:boundingComponentsWithSmallSigma}, \ref{lem:boundingComponentsWithKappa} and \ref{lem:WaldThm1Type}, 
we obtain the following inequality for each $\Psi_{n, \scrk, s}$. 
\begin{eqnarray}
  \limsup_{n \rightarrow \infty}
  \sup_{\theta \in \Psi_{n,\scrm,s}}
  \frac{1}{n}  \sum_{i=1}^{n} 
  \log{f(X_i; \theta)} + \log{r_{n(\theta)}}
 & \leq &
E_{0}[\log{f(x;\theta_0)}]
   - {\lambda_{0}} ,
   \quad a.s.
   \nonumber 
\end{eqnarray}
This leads to (\ref{eq:goal-Psi:thm:penalt-ratio-scale}). 

\paragraph{Setting up constants}

For $\kappa, \lambda$
satisfying~(\ref{eq:kappa-lambda}), let $\kappa_{0}, \lambda_{0}$ be
real constants such that 
\begin{equation}
 0 < 4\kappa_{0} \leq \kappa
  \quad , \quad 
  0 < 4\lambda_{0} \leq \lambda
  \quad , \quad
  \frac{v_{0}}{\kappa_{0}} > \max{\{\sigma_{01}, \ldots, \sigma_{0M}\}}.
  \label{eq:constants-kappa-lambda:penalt-ratio-scale}
\end{equation}
Note that $4\kappa_{0}, 4\lambda_{0}$ also satisfy~(\ref{eq:kappa-lambda}).
Define 
\begin{eqnarray}
 B \equiv \frac{v_{0}}{\kappa_{0}} 
  \label{eq:def:B:thm:penalt-ratio-scale}
\end{eqnarray}
If $\sigma_{m} \geq B$, then the density of the $m$-th component is almost flat 
and makes little contribution to the likelihood. 
In the following argument, we partition the parameter space 
according to this property. 

Because $\{c_n\}$ is decreasing
to zero, by replacing $c_0$ by some $c_n$ if necessary, 
we can assume without loss
of generality that $c_0$ is sufficiently small to satisfy the
following conditions, 
\begin{eqnarray}
 & & (v_0/c_0)^{\tilde{\beta}}
  > e , 
  \nonumber \\
 & & c_{0}  <  \min{\{\sigma_{01}, \ldots, \sigma_{0M}\}} , 
  \nonumber \\
 & & 3M \cdot u_{0} \cdot 2\nu(c_{0}) \cdot \abslr{\log{ 2{\kappa_{0}} }}
   < 
  {\lambda_{0}} ,
   \nonumber 
  \\ 
 & &
 3 \cdot 2M \cdot u_{0} \cdot
  \xi(v_0/c_{0})
  \cdot \log(v_0/c_0)
   < 
  {\lambda_{0}}
  \; , 
   \nonumber 
  \\
  & & {\kappa_{0}}
  < 
  \frac{v_0}{c_0(M + 1)}
  \quad , 
   \label{eq:condition:thm:penalt-ratio-scale}
\end{eqnarray}
where $\tilde{\beta} \equiv (\beta-1)/\beta$ and 
\begin{equation}
 \xi(y) 
  \equiv 
  2\cdot \left(
     \frac{v_{1}}{\kappa_{0}}
    \right)^{\frac{1}{\beta}}
   \cdot 
   \left(
    {v_{0}\cdot {(M + 1)}}
   \right)^{\tilde{\beta}}
   \cdot \left(\frac{1}{y}\right)^{\tilde{\beta}}.
   \label{eq:def:xi:thm:penalt-ratio-scale}
\end{equation}

Take $A_{0} > 0$ sufficiently large such that 
\begin{eqnarray}
  P_{0}( (-\infty, -A_{0}) \cup (A_{0},\infty) )
  \cdot
  \log{
  \left(
   \frac{v_0/c_0 + 3\kappa_{0}}{4\kappa_{0}}
  \right)
  }
  < {\lambda_{0}}. 
  \label{eq:A-zero:condition:thm:penalt-ratio-scale}
\end{eqnarray}
Let
$
 \scra_{0} \equiv (-\infty, -A_{0}) \cup (A_{0},\infty) 
$
and 
$
 A_{n} \equiv A_{0} \cdot n^{\frac{2 + \zeta}{\beta - 1}} 
$
as in Lemma \ref{lem:importantInterval}. 

\paragraph{Partitioning the parameter space}

Partition $\{1,\dots,M\}$ into disjoint subsets 
$\scrm_{\sigma < c_{n}}$, $\scrm_{c_{n} \leq \sigma < c_{0}}$, $\scrm_{\sigma > B}$, $\scrm_{\abs{\mu} > A_{0}}$ and $\scrm_{R}$. 
For any given 
$\scrm_{\sigma < c_{n}}$, $\scrm_{c_{n} \leq \sigma < c_{0}}$, $\scrm_{\sigma > B}$, $\scrm_{\abs{\mu} > A_{0}}$ and $\scrm_{R}$, 
we define a subset of $\Psi_{n}$ by
\begin{eqnarray}
 \Psi_{n, \scrm}
  & \equiv &
  \{
  \theta \in \Psi_{n}
  \mid
  \sigma_{m} < c_{n} , (m \in \scrm_{\sigma < c_{n}}); 
  \nonumber \\ & &
  c_{n} \leq \sigma_{m} < c_{0} , (m \in \scrm_{c_{n} \leq \sigma < c_{0}}); 
  \nonumber \\ & &
  \sigma_{m} > B , (m \in \scrm_{\sigma > B}); 
  \nonumber \\ & & 
  c_{0} \leq \sigma_{m} \leq B , \abs{\mu_{m}} > A_{0} , (m \in \scrm_{\abs{\mu} > A_{0}}); 
  \nonumber \\ & &
  c_{0} \leq \sigma_{m} \leq B , \abs{\mu_{m}} \leq A_{0} , (m \in \scrm_{R}) 
  \} 
  \nonumber 
\end{eqnarray}
The method of partitioning of the parameter space is the same as in Section 4.3.2 of \citet{TT2003-35} 
except for $\scrm_{\sigma < c_{n}}$. 
We will show that the contributions of the components in $\scrm_{\sigma < c_{n}}$ to the likelihood 
are canceled out by the penalty term. 

As above, it suffices to prove that for each choice of disjoint subsets 
$\scrm_{\sigma < c_{n}}$, $\scrm_{c_{n} \leq \sigma < c_{0}}$, $\scrm_{\sigma > B}$, $\scrm_{\abs{\mu} > A_{0}}$ and $\scrm_{R}$ 
\begin{eqnarray}
  \lim_{n \rightarrow \infty}
  \frac{
  \sup_{\theta \in \Psi_{n,\scrm}}
  \prod_{i=1}^{n} f(X_{i};\theta) \cdot r_{n}(\theta)
  }
  { \prod_{i=1}^{n} f(X_{i};\theta_{0}) \cdot r_{n}(\theta_{0})  }
  = 0, \quad a.s. 
  \nonumber 
\end{eqnarray}
We fix 
$\scrm_{\sigma < c_{n}}$, $\scrm_{c_{n} \leq \sigma < c_{0}}$, $\scrm_{\sigma > B}$, $\scrm_{\abs{\mu} > A_{0}}$ and $\scrm_{R}$ 
from now on.

Next we consider coverings of $\bar{\Theta}_{\scrm_{R}}$. 
The following lemma follows immediately from lemma \ref{lem:WaldThm1Type} and compactness of $\bar{\Theta}_{\scrm_{R}}$. 
\begin{lem}
\label{lem:WaldThm1Type:thm:penalt-ratio-scale}
Let $\scrb(\theta,\rho(\theta) )$ denote the open ball
with center $\theta$ and radius $\rho(\theta)$.
Then $\bar{\Theta}_{\scrm_{R}}$ can be covered by 
a finite number of balls 
$\scrb(\theta_{\scrm_{R}}^{(1)}, \rho(\theta_{\scrm_{R}}^{(1)}))
, \ldots, \scrb(\theta_{\scrm_{R}}^{(S)}, \rho(\theta_{\scrm_{R}}^{(S)})) 
$
such that 
\begin{eqnarray}
 E_{0}[\log{\{f_{\scrm_{R}}(x;\theta_{\scrm_{R}}^{(s)},\rho(\theta_{\scrm_{R}}^{(s)})) + \kappa_{0}\}}]
  + \lambda_{0}
  <
  E_{0}[\log{f(x;\theta_0)}]
  \; , \quad (s=1,\ldots,S) \; . 
  \nonumber 
\end{eqnarray}
\end{lem}

\bigskip
Based on lemma \ref{lem:WaldThm1Type:thm:penalt-ratio-scale} we partition $\Psi_{n,\scrm}$.
Recall that we denote by $\theta_{\scrm}$ the subvector of $\theta \in \Theta$ 
consisting of the components in $\scrm$.
Define a subset of $\Psi_{n,\scrm}$ by
\begin{eqnarray}
 \Psi_{n,\scrm, s}
  \equiv 
  \{
  \theta \in \Psi_{n,\scrm}
  \mid 
  \theta_{\scrm_{R}} \in \scrb(\theta_{\scrm_{R}}^{(s)},\rho(\theta_{\scrm_{R}}^{(s)}))
  \}.
  \nonumber 
\end{eqnarray}
Then $\Psi_{n,\scrm}$ is covered by
$\Psi_{n,\scrm,1}, \ldots, \Psi_{n,\scrm,S}$ :
$$
\Psi_{n,\scrm} = \bigcup_{s=1}^S \Psi_{n,\scrm,s} \ .
$$
Again it suffices to prove
that for each choice of 
$\scrm_{\sigma < c_{n}}$, $\scrm_{c_{n} \leq \sigma < c_{0}}$, $\scrm_{\sigma > B}$, $\scrm_{\abs{\mu} > A_{0}}$, $\scrm_{R}$ and $s$
\begin{eqnarray}
   \lim_{n \rightarrow \infty}
   \frac{
   \sup_{\theta \in \Psi_{n,\scrm,s}}
   \prod_{i=1}^{n} f(X_{i};\theta) \cdot r_{n}(\theta)
   }
   { \prod_{i=1}^{n} f(X_{i};\theta_{0}) \cdot r_{n}(\theta_{0}) }
   = 0, \quad a.s.
  \label{eq:goal:PhiLikelihoodRatio}
\end{eqnarray}
We fix 
$\scrm_{\sigma < c_{n}}$, $\scrm_{c_{n} \leq \sigma < c_{0}}$, $\scrm_{\sigma > B}$, $\scrm_{\abs{\mu} > A_{0}}$, $\scrm_{R}$ and $s$
from now on.
By Assumption~\ref{assumption:6},~\ref{assumption:8} and the law of large numbers, 
(\ref{eq:goal:PhiLikelihoodRatio}) is implied by
\begin{eqnarray}
\limsup_{n \rightarrow \infty} \frac{1}{n} \sup_{\theta \in \Psi_{n,\scrm,s}}
   \sum_{i=1}^{n} \log{f(X_{i};\theta)} + \log{r_{n}(\theta)} 
  < 
   E_0[\log f(x;\theta_0)], 
  \quad a.s.
  \label{eq:goal:PhiLogLikelihood}
\end{eqnarray}
Therefore it suffices to prove (\ref{eq:goal:PhiLogLikelihood}). 

\paragraph{Bounding the penalized likelihood by six terms}

The outline of the rest of our proof is as follows. 
First, we bound the likelihood by four terms in Lemma~\ref{lem:boundingLikelihoodFunc-Psi:thm:penalt-ratio-scale}. 
Next, we bound one term of the four terms obtained in Lemma~\ref{lem:boundingLikelihoodFunc-Psi:thm:penalt-ratio-scale} 
by three terms in Lemma~\ref{lem:boundingLikelihoodFunc-Psi2:thm:penalt-ratio-scale}. 
Finally, from Lemma~\ref{lem:boundingLikelihoodFunc-Psi:thm:penalt-ratio-scale} 
and Lemma~\ref{lem:boundingLikelihoodFunc-Psi2:thm:penalt-ratio-scale}, 
we bound the penalized likelihood by six terms 
in Lemma~\ref{lem:boundingLikelihoodFunc-Psi3:thm:penalt-ratio-scale}. 

Define $J_{\sigma < c_{n}}(\theta)$ as 
\begin{eqnarray}
 J_{\sigma < c_{n}}(\theta)
  & \equiv &
  \bigcup_{m \in \scrm_{\sigma < c_{n}}}[\mu_{m}-\nu(\sigma_{m}), \; \mu_{m}+\nu(\sigma_{m})]. 
 \label{eq:definitionOfJSigmaLeqCn:thm:penalt-ratio-scale} 
\end{eqnarray}
Let $\scrm_{\sigma \geq c_{n}}=\{1,\dots,M\}\setminus\scrm_{\sigma < c_{n}}$. 
Then the following lemma holds. 
\begin{lem} For $\theta \in \Psi_{n,\scrm,s}$, 
 \label{lem:boundingLikelihoodFunc-Psi:thm:penalt-ratio-scale}
 \begin{eqnarray}
 \lefteqn{
  \sum_{i=1}^{n}\log{f(X_{i};\theta)}
  } & & \nonumber \\
 & \leq & 
   \sum_{i=1}^{n}
   \log{
   \left\{
    \sum_{m \in \scrm_{\sigma \geq c_{n}}}\alpha_{m}f_{m}(X_{i};\theta_{m})
    + 
    \kappa_{0}
   \right\}
   }
  + 
   R_{n}(J_{\sigma < c_{n}}(\theta)) \cdot \log{\frac{v_{0}}{\kappa_{0}}}
  \nonumber \\ 
 & & + 
   R_{n}(J_{\sigma < c_{n}}(\theta)) \cdot \log{\frac{1}{\sigma_{(1)}}}. 
  \label{eq:boundingLikelihoodFunc-Psi:lem:boundingLikelihoodFunc-Psi:thm:penalt-ratio-scale}
 \end{eqnarray}
\end{lem}
\Proof
For $\theta \in \Psi_{n,\scrm,s} \subset \Psi_{n} \subset \Theta_{c_n}^{C}$, 
from Lemma~\ref{lem:boundingComponentsWithKappa}, the following inequalities hold. 
\begin{eqnarray}
 \lefteqn{
 \sum_{i=1}^{n}\log{f(X_{i};\theta)} 
 =
   \sum_{X_{i} \in J_{\sigma < c_{n}}(\theta)} \log{f(X_{i};\theta)} 
  + 
   \sum_{X_{i} \in \real\setminus J_{\sigma < c_{n}}(\theta)} \log{f(X_{i};\theta)}   
   } \qquad \qquad \qquad \qquad & & \nonumber \\
 & \leq & 
   R_{n}(J_{\sigma < c_{n}}(\theta)) \cdot \log{\left\{\max_{1 \leq m \leq M}\left(\frac{v_{0}}{\sigma_{m}}\right)\right\}}
  \nonumber \\ 
   & & + 
   \sum_{X_{i} \in \real\setminus J_{\sigma < c_{n}}(\theta)} 
   \log{
   \left\{
    \sum_{m \in \scrm_{\sigma \geq c_{n}}}\alpha_{m}f_{m}(X_{i};\theta_{m})
    + 
    \kappa_{0}
   \right\}
   } 
  \nonumber \\ 
 & \leq & 
  \sum_{i=1}^{n}
  \log{
  \left\{
   \sum_{m \in \scrm_{\sigma \geq c_{n}}}\alpha_{m}f_{m}(X_{i};\theta_{m})
   + 
   \kappa_{0}
  \right\}
  }
  - 
  \sum_{X_{i} \in J_{\sigma < c_{n}}(\theta)} \log{\kappa_{0}}
  \nonumber \\ 
 & & + 
   R_{n}(J_{\sigma < c_{n}}(\theta)) \cdot \log{\frac{v_{0}}{\sigma_{(1)}}}
  \nonumber \\
 & \leq & 
  \sum_{i=1}^{n}
  \log{
  \left\{
   \sum_{m \in \scrm_{\sigma \geq c_{n}}}\alpha_{m}f_{m}(X_{i};\theta_{m})
   + 
   \kappa_{0}
  \right\}
  }
  + 
   R_{n}(J_{\sigma < c_{n}}(\theta)) \cdot \log{\frac{v_{0}}{\kappa_{0}}}
  \nonumber \\ 
 & & + 
   R_{n}(J_{\sigma < c_{n}}(\theta)) \cdot \log{\frac{1}{\sigma_{(1)}}}. 
   \nonumber 
\end{eqnarray}
\qed

Define $J_{c_{n} \leq \sigma < c_{0}}(\theta)$ as 
\begin{eqnarray}
 J_{c_{n} \leq \sigma < c_{0}}(\theta)
  & \equiv &
  \bigcup_{m \in \scrm_{c_{n} \leq \sigma < c_{0}}}[\mu_{m}-\nu(\sigma_{m}), \; \mu_{m}+\nu(\sigma_{m})]. 
  \nonumber 
\end{eqnarray}
For the first term of (\ref{eq:boundingLikelihoodFunc-Psi:lem:boundingLikelihoodFunc-Psi:thm:penalt-ratio-scale}), 
the following lemma holds. 
\begin{lem}
  \label{lem:boundingLikelihoodFunc-Psi2:thm:penalt-ratio-scale}
 The following inequality holds for $\theta \in \Psi_{n,\scrm,s}$. 
 \begin{eqnarray}
  \sum_{i=1}^{n}
   \log{
   \left\{
    \sum_{m \in \scrm_{\sigma \geq c_{n}}}\alpha_{m}f_{m}(X_{i};\theta_{m})
    + 
    \kappa_{0}
   \right\}
   }
 & \leq & 
  \sum_{i=1}^{n}
  \log{
  \left\{
   f_{\scrm_{R}}(X_{i};\theta_{\scrm_{R}},\rho(\theta_{\scrm_{R}})) + 4{\kappa_{0}}
  \right\}
  }
  \nonumber \\
 & & + \; 
  R_{n}(\scra_{0}) \cdot
  \log{
  \left(
   \frac{v_0/c_0 + 3\kappa_{0}}{4\kappa_{0}}
  \right)
  } 
  \nonumber \\
 & & + \; 
  R_n(J_{c_{n} \leq \sigma < c_{0}}(\theta))\cdot 
  (
  -\log{2\kappa_{0}} 
  )
  \nonumber \\
 & & + \; 
  \sum_{X_{i} \in J_{c_{n} \leq \sigma < c_{0}}(\theta)}
  \log{\{f(X_{i};\theta) + \kappa_{0}\}}
  \; . 
  \nonumber \\
  \label{eq:boundingLikelihoodFunc-Psi2:lem:boundingLikelihoodFunc-Psi2:thm:penalt-ratio-scale}
 \end{eqnarray}
\end{lem}
\Proof
Let $\scrm_{{\sigma} \geq c_{0}} \equiv \{1,\dots,M\}\backslash \{ \scrm_{\sigma < c_{n}} \cup \scrm_{c_{n} \leq \sigma < c_{0}} \}$ 
and 
$\scrm_{c_{n} \leq \sigma \leq B} \equiv \{1,\dots,M\}\backslash \{ \scrm_{\sigma < c_{n}} \cup \scrm_{c_{n} \leq \sigma < c_{0}} \cup \scrm_{{\sigma} > B}\}$.
For $x \not\in J_{c_{n} \leq \sigma < c_{0}}(\theta)$, 
$f(x;\theta)\leq f_{\scrm_{{\sigma} > c_{0}}}
(x;\theta_{\scrm_{\sigma \geq c_{0}}}) 
+ \kappa_{0}$ 
holds.
Therefore
\begin{eqnarray}
 \sum_{i=1}^{n}
  \log{\left\{ 
	f(X_{i};\theta) + \kappa_{0}
	\right\}}
  & \leq & 
   \sum_{X_i \in J_{c_{n} \leq \sigma < c_{0}}(\theta)} 
  \log{\left\{ 
	f(X_{i};\theta) + \kappa_{0}
       \right\}}
  \nonumber \\ 
 & &+ 
   \sum_{X_i \not \in J_{c_{n} \leq \sigma < c_{0}}(\theta)} 
  \log{\left\{ 
       f_{\scrm_{\sigma \geq c_{0}}}(x;\theta_{\scrm_{\sigma \geq c_{0}}}) + 2\kappa_{0}
      \right\}}
  \nonumber \\
 &= &
    \sum_{i=1}^{n}
  \log{
  \left\{
   f_{\scrm_{\sigma \geq c_{0}}}(x;\theta_{\scrm_{\sigma \geq c_{0}}}) + 2\kappa_{0}
  \right\}
  }
  \nonumber \\
 && +
  \sum_{X_{i} \in J_{c_{n} \leq \sigma < c_{0}}(\theta)}
  \left[
   \log{\left\{ 
	 f(X_{i};\theta) + \kappa_{0}
	\right\}}
   -
   \log{
   \left\{
    f_{\scrm_{\sigma \geq c_{0}}}(x;\theta_{\scrm_{\sigma \geq c_{0}}}) + 2\kappa_{0}
   \right\}}
  \right]
  \nonumber \\
\label{eq:lemma9proof1}
\end{eqnarray}
Consider the second term on the right-hand side. We have
\begin{align*}
&  \sum_{X_{i} \in J_{c_{n} \leq \sigma < c_{0}}(\theta)}
  \left[
   \log{\left\{ 
	 f(X_{i};\theta) + \kappa_{0}
	\right\}}
   -
   \log{
   \left\{
    f_{\scrm_{\sigma \geq c_{0}}}(x;\theta_{\scrm_{\sigma \geq c_{0}}}) + 2\kappa_{0}
   \right\}}
  \right]\\
& \qquad\qquad \le 
   \sum_{X_{i} \in J_{c_{n} \leq \sigma < c_{0}}(\theta)}
   \log{\left\{ 
	 f(X_{i};\theta) + \kappa_{0}
	\right\}}
   -
   R_n(J_{c_{n} \leq \sigma < c_{0}}(\theta)) \cdot \log{2\kappa_{0}}
   \; . 
\end{align*}
This takes care of the third and the fourth term of (\ref{eq:boundingLikelihoodFunc-Psi2:lem:boundingLikelihoodFunc-Psi2:thm:penalt-ratio-scale}).
Now consider the first term on the right-hand side of
(\ref{eq:lemma9proof1}).
Because of (\ref{eq:def:B:thm:penalt-ratio-scale}), we obtain
\begin{eqnarray}
    \sum_{i=1}^{n}
  \log{
  \left\{
   f_{\scrm_{\sigma \geq c_{0}}}(X_{i};\theta_{\scrm_{\sigma \geq c_{0}}}) + 2\kappa_{0}
  \right\}
  }
   &\leq &
  \sum_{i=1}^{n}
  \log{
  \left\{
     f_{\scrm_{c_{0} \leq \sigma \leq B}}(X_{i};\theta_{\scrm_{c_{0} \leq \sigma \leq B}}) + 3\kappa_{0}
  \right\}
  }. 
  \nonumber 
\end{eqnarray}
Note that $\scra_{0} = \{x \in \real \mid \abs{x} > A_{0}\}$ and 
$\scrm_{R} = \scrm_{c_{0} \leq \sigma \leq B} \setminus \scrm_{\abs{\mu} > A_{0}}$. 
For $x\notin  \scra_{0}$, we have 
$$
 f_{\scrm_{\abs{\mu} > A_{0}}}(x;\theta_{\scrm_{\abs{\mu} > A_{0}}})
  \leq
 {\kappa_{0}}.
$$
Therefore we obtain 
\begin{eqnarray}
 \lefteqn{
 \sum_{i=1}^{n}
  \log{
  \left\{
     f_{\scrm_{c_{0} \leq \sigma \leq B}}(X_{i};\theta_{\scrm_{c_{0} \leq \sigma \leq B}})
   + 3\kappa_{0}
  \right\}
  } }
  \nonumber \\
  & = & 
  \sum_{X_{i} \notin \scra_{0}}
  \log{
  \left\{
   f_{\scrm_{c_{0} \leq \sigma \leq B}}(X_{i};\theta_{\scrm_{c_{0} \leq \sigma \leq B}})
   + 3\kappa_{0}
  \right\}
  }
  + 
  \sum_{X_{i} \in \scra_{0}}
  \log{
  \left\{
   f_{\scrm_{c_{0} \leq \sigma \leq B}}(X_{i};\theta_{\scrm_{c_{0} \leq \sigma \leq B}})
   + 3\kappa_{0}
  \right\}
  }
  \nonumber \\
 & \leq &
    \sum_{X_{i} \notin \scra_{0}}
  \log{
  \left\{
   f_{\scrm_{R}}(X_{i};\theta_{\scrm_{R}})
   + 4\kappa_{0}
  \right\}
  }
  + 
  \sum_{X_{i} \in \scra_{0}}
  \log{
  \left\{
   f_{\scrm_{c_{0} \leq \sigma \leq B}}(X_{i};\theta_{\scrm_{c_{0} \leq \sigma \leq B}})
   + 3\kappa_{0}
  \right\}
  }
  \nonumber \\
 & = & 
  \sum_{i=1}^{n}
  \log{
  \left\{
   f_{\scrm_{R}}(X_{i};\theta_{\scrm_{R}})
   + 4\kappa_{0}
  \right\}
  }
  \nonumber \\
 & & \qquad \qquad 
  + 
  \sum_{X_{i} \in \scra_{0}}
  \left[
   \log{
   \left\{
    f_{\scrm_{c_{0} \leq \sigma \leq B}}(X_{i};\theta_{\scrm_{c_{0} \leq \sigma \leq B}})
    + 3\kappa_{0}
   \right\}
   }
  - 
   \log{
   \left\{
   f_{\scrm_{R}}(X_{i};\theta_{\scrm_{R}})
    + 4\kappa_{0}
   \right\}
   }
  \right]
  \nonumber  \\ 
  \label{eq:totyu}
\end{eqnarray}
Note that $f_{\scrm_{c_{0} \leq \sigma \leq B}}(x;\theta_{\scrm_{c_{0} \leq \sigma \leq B}}) \leq v_{0}/c_{0}$
from lemma~\ref{lem:boundingComponentsWithSmallSigma}. Therefore 
\begin{eqnarray}
 \lefteqn{
 \text{The r.h.s of (\ref{eq:totyu})} 
 } & & \nonumber \\
 &\leq &
  \sum_{i=1}^{n}
  \log{
  \left\{
   f_{\scrm_{R}}(X_{i};\theta_{\scrm_{R}})
   + 4\kappa_{0}
  \right\}
  }
  + 
  \sum_{X_{i} \in \scra_{0}}
  \left[
   \log{
   \left\{
    v_0/c_0
    + 3\kappa_{0}
   \right\}
   }
   - 
   \log{
   4\kappa_{0}
   }
  \right]
  \nonumber \\
 &\leq &
  \sum_{i=1}^{n}
  \log{
  \left\{
     f_{\scrm_{R}}(X_{i};\theta_{\scrm_{R}},\rho(\theta_{\scrm_{R}})) + 4{\kappa_{0}}
  \right\}
  }
  + 
   R_{n}(\scra_{0}) \cdot 
  \log{
  \left(
   \frac{v_0/c_0 + 3{\kappa_{0}}}{4{\kappa_{0}}}
  \right)
  }.
  \nonumber 
\end{eqnarray}
This takes care of the first and the second term of (\ref{eq:boundingLikelihoodFunc-Psi2:lem:boundingLikelihoodFunc-Psi2:thm:penalt-ratio-scale}). 
\qed

By Lemma~\ref{lem:boundingLikelihoodFunc-Psi:thm:penalt-ratio-scale} and \ref{lem:boundingLikelihoodFunc-Psi2:thm:penalt-ratio-scale}, 
the log likelihood function is bounded above as the following lemma. 
\begin{lem}
 \label{lem:boundingLikelihoodFunc-Psi3:thm:penalt-ratio-scale}
 For $\theta \in \Psi_{n,\scrm,s}$, 
  \begin{eqnarray}
  \lefteqn{
  \sum_{i=1}^{n}\log{f(X_{i};\theta)} + \log{r_{n}(\theta)}
  } & & 
  \nonumber \\
  & \leq  & 
  \sum_{i=1}^{n}
  \log{
  \left\{
   f_{\scrm_{R}}(X_{i};\theta_{\scrm_{R}},\rho(\theta_{\scrm_{R}})) + 4{\kappa_{0}}
  \right\}
  }
  \nonumber \\
 & & +
  R_{n}(\scra_{0}) \cdot
  \log{
  \left(
   \frac{v_0/c_0 + 3\kappa_{0}}{4\kappa_{0}}
  \right)
  } 
  \nonumber \\
 & & +
  R_n(J_{c_{n} \leq \sigma < c_{0}}(\theta))\cdot 
  (
  -\log{2\kappa_{0}} 
  )
  \nonumber \\
 & & +
  \sum_{X_{i} \in J_{c_{n} \leq \sigma < c_{0}}(\theta)}
  \log{\{f(X_{i};\theta) + \kappa_{0}\}}
  \nonumber \\
  & &   + 
   R_{n}(J_{\sigma < c_{n}}(\theta)) \cdot \log{\frac{v_{0}}{\kappa_{0}}}
  \nonumber \\ 
 & & + 
  \left\{
   R_{n}(J_{\sigma < c_{n}}(\theta)) \cdot \log{\frac{1}{\sigma_{(1)}}}
   + 
   \log{r_{n}(\theta)} 
  \right\}. 
  \label{eq:boundingLikelihoodFunc-Psi3:lem:boundingLikelihoodFunc-Psi3:thm:penalt-ratio-scale}
 \end{eqnarray}
\end{lem}

We bound the six terms of (\ref{eq:boundingLikelihoodFunc-Psi3:lem:boundingLikelihoodFunc-Psi3:thm:penalt-ratio-scale}) 
in the following paragraphs. 

\paragraph{The first term}
We begin by bounding the first term of (\ref{eq:boundingLikelihoodFunc-Psi3:lem:boundingLikelihoodFunc-Psi3:thm:penalt-ratio-scale}). 
By lemma~\ref{lem:WaldThm1Type} and the strong law of large
numbers, we have
\begin{eqnarray}
\limsup_{n\rightarrow\infty}  \frac{1}{n}\sum_{i=1}^{n}
  \log{
  \left\{
     f_{\scrm_{R}}(X_{i};\theta_{\scrm_{R}},\rho(\theta_{\scrm_{R}})) + 4{\kappa_{0}}
  \right\}
  }
 < 
 E_0[\log f(x;\theta_0)] - 4\lambda_{0} , \quad a.s.
  \label{eq:BoundingTheFirstTerm:thm:penalt-ratio-scale}
\end{eqnarray}

\paragraph{The second term}
By (\ref{eq:A-zero:condition:thm:penalt-ratio-scale}) and the strong law of large
numbers, we have
\begin{eqnarray}
  \limsup_{n\rightarrow\infty} 
  \frac{1}{n}R_{n}(\scra_{0}) \cdot
  \log{
  \left(
   \frac{v_0/c_0 + 3\kappa_{0}}{4\kappa_{0}}
  \right)
  } 
  < 
  {\lambda_{0}}
  , \quad a.s.
  \label{eq:BoundingTheSecondTerm:thm:penalt-ratio-scale}
\end{eqnarray}

\paragraph{The third term and the fourth term}
The third term and the fourth term of (\ref{eq:boundingLikelihoodFunc-Psi3:lem:boundingLikelihoodFunc-Psi3:thm:penalt-ratio-scale}) 
can be bounded from above as follows: 
 \begin{eqnarray}
  & &
   \limsup_{n \rightarrow \infty}
   \sup_{\theta \in \Psi_{n,\scrm,s}}
   \frac{1}{n} R_{n}(J_{c_{n} \leq \sigma < c_{0}}(\theta))\cdot \abs{\log{2\kappa_{0}}}
   \leq 3M \cdot u_{0} \cdot 2\nu(c_{0}) \cdot \abs{\log{2\kappa_{0}}} < \lambda_{0}, 
   \quad a.s.
   \nonumber \\
   \label{eq:BoundingTheThirdTerm:thm:penalt-ratio-scale}
\end{eqnarray}
\begin{eqnarray} 
   \limsup_{n \rightarrow \infty}
   \sup_{\theta \in \Psi_{n,\scrk,s}}
   \frac{1}{n}  \sum_{X_{i} \in J_{c_{n} \leq \sigma < c_{0}}(\theta)}
   \log{
   \left\{
    f(X_i; \theta) + \kappa_{0}
   \right\}
   }
   \leq
   {\lambda_{0}} ,
   \quad a.s.
   \nonumber \\
  \label{eq:BoundingTheFourthTerm:thm:penalt-ratio-scale}
 \end{eqnarray}
The proofs of the above inequalities are similar to the proofs of section 4.3.4 and 4.3.5 
in the longer version of \citet{TT2003-35}, and are omitted.

\paragraph{The fifth term and the sixth term}
We now state the following lemma 
in order to bound the fifth term and sixth term of (\ref{eq:boundingLikelihoodFunc-Psi3:lem:boundingLikelihoodFunc-Psi3:thm:penalt-ratio-scale}). 
\begin{lem}
 \label{lem:numberOfObservationsWithinShortInterval-2:thm:penalt-ratio-scale}
 \begin{equation}
  \prob\left(
	\sup_{\theta \in \Psi_{n,\scrm,s}}
	R_{n}(J_{\sigma < c_{n}}(\theta)) > M \quad \io
       \right) = 0 
  \label{eq:lem:numberOfObservationsWithinShortInterval-2:thm:penalt-ratio-scale}
 \end{equation}
\end{lem}
\Proof
Let $w_{n} = \nu(c_{n}) = \nu(\exp(n^{-\cnd}))$. 
Then (\ref{eq:assumption:lem:nondense}), the assumption of Lemma~\ref{lem:nondense}, is satisfied. 
From (\ref{eq:definitionOfJSigmaLeqCn:thm:penalt-ratio-scale}), we have 
\begin{eqnarray}
  \sup_{\theta \in \Psi_{{n},\scrm,s}}
  R_{n}(J_{\sigma < c_{n}}(\theta)) > M 
  \quad & \Rightarrow & \quad 
  \max_{\mu \in \real}R_{n}([\mu - w_{n},\, \mu + w_{n}] )  > 1 
  \nonumber 
\end{eqnarray}
Therefore, by Lemma~\ref{lem:nondense}, we obtain (\ref{eq:lem:numberOfObservationsWithinShortInterval-2:thm:penalt-ratio-scale}). 
\qed

By Lemma~\ref{lem:numberOfObservationsWithinShortInterval-2:thm:penalt-ratio-scale} 
and the same argument in Section~\ref{sec:proof:eq:goal-Phi:thm:penalt-ratio-scale}, 
we ignore the event \\
$R_{n}(J_{\sigma < c_{n}}(\theta)) > M$. 
Then we have for $\theta \in \Psi_{n, \scrm, s}$ uniformly
\begin{eqnarray}
   R_{n}(J_{\sigma < c_{n}}(\theta)) \cdot \log{\frac{v_{0}}{\kappa_{0}}}
    + 
  \left\{
   R_{n}(J_{\sigma < c_{n}}(\theta)) \cdot \log{\frac{1}{\sigma_{(1)}}}
   + 
   \log{r_{n}(\theta)} 
  \right\}
  \nonumber \\
 \leq 
   M \cdot \log{\frac{v_{0}}{\kappa_{0}}}
  + 
   \log{\frac{r_{n}(\theta)}{(\sigma_{(1)})^{M}}}
  \qquad , \qquad  
  a.s. 
  \nonumber 
\end{eqnarray}
From~(\ref{eq:definitionOfPsi:thm:penalt-ratio-scale}), we obtain for $\theta \in \Psi_{n, \scrm, s}$
\begin{eqnarray}
 R_{n}(J_{\sigma < c_{n}}(\theta)) \cdot \log{\frac{v_{0}}{\kappa_{0}}}
 + 
 \left\{
  R_{n}(J_{\sigma < c_{n}}(\theta)) \cdot \log{\frac{1}{\sigma_{(1)}}}
  + 
  \log{r_{n}(\theta)} 
 \right\}
 \nonumber \\
 \leq
   M \cdot \log{\frac{v_{0}}{\kappa_{0}}}
  \qquad , \qquad a.s. 
  \nonumber 
\end{eqnarray}

Because the right hand side of the above inequality is constant, we have 
\begin{eqnarray}
 \limsup_{n\rightarrow \infty}
  \frac{1}{n} \cdot \sup_{\theta \in \Psi_{n, \scrm ,s}}
  \left[
    R_{n}(J_{\sigma < c_{n}}(\theta)) \cdot \log{\frac{v_{0}}{\kappa_{0}}}
   + 
   \left\{
    R_{n}(J_{\sigma < c_{n}}(\theta)) \cdot \log{\frac{1}{\sigma_{(1)}}}
   + 
   \log{r_{n}(\theta)} 
   \right\}
   \right]
  \nonumber \\
  \leq 0 \quad a.s. 
  \nonumber \\
 \label{eq:BoundingTheFifthAndSixthTerm:thm:penalt-ratio-scale}
\end{eqnarray}

\paragraph{The end of the proof}

Combining 
(\ref{eq:boundingLikelihoodFunc-Psi3:lem:boundingLikelihoodFunc-Psi3:thm:penalt-ratio-scale}), 
(\ref{eq:BoundingTheFirstTerm:thm:penalt-ratio-scale}), 
(\ref{eq:BoundingTheSecondTerm:thm:penalt-ratio-scale}), 
(\ref{eq:BoundingTheThirdTerm:thm:penalt-ratio-scale}), 
(\ref{eq:BoundingTheFourthTerm:thm:penalt-ratio-scale}), 
(\ref{eq:BoundingTheFifthAndSixthTerm:thm:penalt-ratio-scale}),
we obtain 
\begin{eqnarray}
  \limsup_{n \rightarrow \infty}
  \sup_{\theta \in \Psi_{n,\scrm,s}}
  \frac{1}{n}  \sum_{i=1}^{n} 
  \log{f(X_i; \theta)} + \log{r_{n(\theta)}}
 & \leq &
E_{0}[\log{f(x;\theta_0)}]
   - {\lambda_{0}} ,
   \quad a.s.
   \nonumber 
\end{eqnarray}
Therefore we obtain (\ref{eq:goal:PhiLogLikelihood}). 

This completes the proof of Theorem~\ref{thm:penalt-ratio-scale}.

\subsection{Proof of Theorem~\ref{thm:penalt-scale-param}}
\label{sec:proof-thm:penalt-scale-param}

The outline of the proof of Theorem~\ref{thm:penalt-scale-param} is 
similar to the proof of Theorem~\ref{thm:penalt-ratio-scale}. 

\paragraph{Partitioning the parameter space}
Let $\cnd$ be a constant defined by Assumption~\ref{assumption:11}. 
Define $c_{n}=c_{0}\cdot\exp(-n^{\cnd})$ and $\Theta_{c_{n}}=\{\theta \in \Theta \mid \sigma_{(1)} \geq c_{n}\}$. 
The parameter space $\Theta$ is divided into two sets.  
\begin{equation}
 \Theta = \Theta_{c_{n}} \cup \Theta_{c_{n}}^{C}. 
  \nonumber 
\end{equation}

Because the asymptotic behavior is not affected by the penalty term, 
the penalized maximum likelihood estimator over $\Theta_{c_{n}}$ is strongly consistent  by Lemma~\ref{thm:constr-scale-param}. 
Therefore, it suffices to prove the following equation.  
\begin{equation}
 \lim_{n \rightarrow \infty}
 \frac{\sup_{\theta \in \Theta_{c_{n}}^{C}}
  \left\{\prod_{i=1}^{n} f(X_{i};\theta) \right\} \cdot s_{n}(\theta)}
  {\left\{ \prod_{i=1}^{n} f(X_{i};\theta_{0}) \right\} \cdot s_{n}(\theta_{0}) } = 0, \quad a.s. 
  \label{eq:lem:goal:thm:penalt-scale-param}
\end{equation}

\paragraph{Setting up constants}

We set up some constants as in section \ref{sec:goal-Psi:thm:penalt-ratio-scale}.

Let $\kappa_{0}, \lambda_{0}$ be real constants such that (\ref{eq:constants-kappa-lambda:penalt-ratio-scale}) holds. 
We can assume without loss
of generality that $c_0$ is sufficiently small to satisfy the equations (\ref{eq:condition:thm:penalt-ratio-scale}). 
Take $A_{0} > 0$ sufficiently large such that (\ref{eq:A-zero:condition:thm:penalt-ratio-scale}) holds. 
Let
$
 \scra_{0} \equiv (-\infty, -A_{0}) \cup (A_{0},\infty) 
$
and 
$
 A_{n} \equiv A_{0} \cdot n^{\frac{2 + \zeta}{\beta - 1}} 
$
as in lemma \ref{lem:importantInterval}. 
Remember that $\tilde{\beta} = (\beta-1)/\beta$, and 
$B$ and $\xi$ are defined in (\ref{eq:def:B:thm:penalt-ratio-scale}) and (\ref{eq:def:xi:thm:penalt-ratio-scale}) respectively. 

\paragraph{Partitioning the parameter space}

Partition \{1, \dots, M\} into disjoint subsets 
$\scrm_{\sigma < c_{n}}$, $\scrm_{c_{n} \leq \sigma < c_{0}}$, $\scrm_{\sigma > B}$, $\scrm_{\abs{\mu} > A_{0}}$ and $\scrm_{R}$. 
For any given 
$\scrm_{\sigma < c_{n}}$, $\scrm_{c_{n} \leq \sigma < c_{0}}$, $\scrm_{\sigma > B}$, $\scrm_{\abs{\mu} > A_{0}}$ and $\scrm_{R}$, 
we define a subset of $\Theta_{c_{n}}^{C}$ by
\begin{eqnarray}
 \Theta_{c_{n}, \scrm}^{C}
  & \equiv &
  \{
  \theta \in \Theta_{c_{n}}^{C}
  \mid
  \sigma_{m} < c_{n} , (m \in \scrm_{\sigma < c_{n}}); 
  \nonumber \\ & &
  c_{n} \leq \sigma_{m} < c_{0} , (m \in \scrm_{c_{n} \leq \sigma < c_{0}}); 
  \nonumber \\ & &
  \sigma_{m} > B , (m \in \scrm_{\sigma > B}); 
  \nonumber \\ & & 
  c_{0} \leq \sigma_{m} \leq B , \abs{\mu_{m}} > A_{0} , (m \in \scrm_{\abs{\mu} > A_{0}}); 
  \nonumber \\ & &
  c_{0} \leq \sigma_{m} \leq B , \abs{\mu_{m}} \leq A_{0} , (m \in \scrm_{R}) 
  \} 
  \nonumber 
\end{eqnarray}
As above, it suffices to prove that for each choice of disjoint subsets 
$\scrm_{\sigma < c_{n}}$, $\scrm_{c_{n} \leq \sigma < c_{0}}$, $\scrm_{\sigma > B}$, $\scrm_{\abs{\mu} > A_{0}}$ and $\scrm_{R}$ 
\begin{eqnarray}
  \lim_{n \rightarrow \infty}
  \frac{
  \sup_{\theta \in \Theta_{c_{n}, \scrm}^{C}}
  \{\prod_{i=1}^{n} f(X_{i};\theta)\}\cdot s_{n}(\theta)
  }
  { \{\prod_{i=1}^{n} f(X_{i};\theta_{0})\}\cdot s_{n}(\theta_{0}) }
  = 0, \quad a.s. 
  \nonumber 
\end{eqnarray}
We fix 
$\scrm_{\sigma < c_{n}}$, $\scrm_{c_{n} \leq \sigma < c_{0}}$, $\scrm_{\sigma > B}$, $\scrm_{\abs{\mu} > A_{0}}$ and $\scrm_{R}$ 
from now on.

Next we consider coverings of $\bar{\Theta}_{\scrm_{R}}$. 
The following lemma follows immediately from lemma \ref{lem:WaldThm1Type} and compactness of $\bar{\Theta}_{\scrm_{R}}$. 
\begin{lem}
\label{lem:WaldThm1Type:thm:penalt-scale-param}
Let $\scrb(\theta,\rho(\theta) )$ denote the open ball
with center $\theta$ and radius $\rho(\theta)$.
Then $\bar{\Theta}_{\scrm_{R}}$ can be covered by 
a finite number of balls 
$\scrb(\theta_{\scrm_{R}}^{(1)}, \rho(\theta_{\scrm_{R}}^{(1)}))
, \ldots, \scrb(\theta_{\scrm_{R}}^{(S)}, \rho(\theta_{\scrm_{R}}^{(S)}))
$
such that 
\begin{eqnarray}
 E_{0}[\log{\{f_{\scrm_{R}}(x;\theta_{\scrm_{R}}^{(s)},\rho(\theta_{\scrm_{R}}^{(s)})) + \kappa_{0}\}}]
  + \lambda_{0}
  <
  E_{0}[\log{f(x;\theta_0)}]
  \; , \quad (s=1,\ldots,S) \; . 
  \nonumber 
\end{eqnarray}
\end{lem}
Based on lemma \ref{lem:WaldThm1Type:thm:penalt-scale-param}, we partition $\Theta_{c_{n}, \scrm}^{C}$.
Define a subset of $\Theta_{c_{n}, \scrm}^{C}$ by
\begin{eqnarray}
 \Theta_{c_{n},\scrm, s}^{C}
  \equiv 
  \{
  \theta \in \Theta_{c_{n}, \scrm}^{C}
  \mid 
  \theta_{\scrm_{R}} \in \scrb(\theta_{\scrm_{R}}^{(s)},\rho(\theta_{\scrm_{R}}^{(s)}))
  \}.
  \nonumber 
\end{eqnarray}
Then $\Theta_{c_{n}, \scrm}^{C}$ is covered by
$\Theta_{c_{n},\scrm,1}^{C}, \ldots, \Theta_{c_{n},\scrm,S}^{C}$ :
$$
\Theta_{c_{n}, \scrm}^{C} = \bigcup_{s=1}^S \Theta_{c_{n},\scrm,s}^{C} \   .
$$
Again it suffices to prove
that for each choice of 
$\scrm_{\sigma < c_{n}}$, $\scrm_{c_{n} \leq \sigma < c_{0}}$, $\scrm_{\sigma > B}$, $\scrm_{\abs{\mu} > A_{0}}$, $\scrm_{R}$ and $s$
\begin{eqnarray}
   \lim_{n \rightarrow \infty}
   \frac{
   \sup_{\theta \in \Theta_{c_{n},\scrm,s}^{C}}
   \{\prod_{i=1}^{n} f(X_{i};\theta)\}\cdot s_{n}(\theta)
   }
   { \{\prod_{i=1}^{n} f(X_{i};\theta_{0})\}\cdot s_{n}(\theta_{0}) }
   = 0, \quad a.s.
  \nonumber 
\end{eqnarray}
We fix 
$\scrm_{\sigma < c_{n}}$, $\scrm_{c_{n} \leq \sigma < c_{0}}$, $\scrm_{\sigma > B}$, $\scrm_{\abs{\mu} > A_{0}}$, $\scrm_{R}$ and $s$
from now on.

By Assumption~\ref{assumption:10},~\ref{assumption:12} and law of large numbers, 
(\ref{eq:lem:goal:thm:penalt-scale-param}) is implied by 
\begin{equation}
 \limsup_{n\rightarrow \infty} \frac{1}{n} \cdot \sup_{\theta \in \Theta_{c_{n},\scrm,s}^{C}}
  \left\{ \sum_{i=1}^{n}\log{f(X_{i};\theta)} + \log{s_{n}(\theta)} \right\}
  < E_{0}[\log{f(x;\theta_{0})}] \quad a.s. 
  \label{eq:goal2:thm:penalt-scale-param}
\end{equation}
We prove (\ref{eq:goal2:thm:penalt-scale-param}) in the following paragraphs. 

\paragraph{Bounding the penalized likelihood function by six terms}

Define $J_{\sigma < c_{n}}(\theta)$ and $J_{c_{n} \leq \sigma < c_{0}}$ as 
\begin{eqnarray}
 J_{\sigma < c_{n}}(\theta) 
  & \equiv &
  \bigcup_{m \in \scrm_{\sigma < c_{n}}}[\mu_{m}-\nu(\sigma_{m}), \; \mu_{m}+\nu(\sigma_{m})].  
 \nonumber \\
 J_{c_{n} \leq \sigma < c_{0}}(\theta) 
  & \equiv &
  \bigcup_{m \in \scrm_{c_{n} \leq \sigma < c_{0}}}[\mu_{m}-\nu(\sigma_{m}), \; \mu_{m}+\nu(\sigma_{m})]. 
 \label{eq:def:J:thm:penalt-scale-param}
\end{eqnarray}
The following lemma can be proved 
by a method similar to the proof of Lemma~\ref{lem:boundingLikelihoodFunc-Psi3:thm:penalt-ratio-scale}. 
\begin{lem}
 \label{lem:boundingLikelihoodFunc:thm:penalt-scale-param}
 For $\theta \in \Theta_{c_{n},\scrm,s}^{C}$, 
  \begin{eqnarray}
  \lefteqn{
  \sum_{i=1}^{n}\log{f(X_{i};\theta)} + \log{s_{n}(\theta)}
  } & & 
  \nonumber \\
  & \leq  & 
  \sum_{i=1}^{n}
  \log{
  \left\{
   f_{\scrm_{R}}(X_{i};\theta_{\scrm_{R}},\rho(\theta_{\scrm_{R}})) + 4{\kappa_{0}}
  \right\}
  }
  \nonumber \\
 & & +
  R_{n}(\scra_{0}) \cdot
  \log{
  \left(
   \frac{v_0/c_0 + 3\kappa_{0}}{4\kappa_{0}}
  \right)
  } 
  \nonumber \\
 & & +
  R_{n}(J_{c_{n} \leq \sigma < c_{0}}(\theta))\cdot 
  (
  -\log{2\kappa_{0}} 
  )
  \nonumber \\
 & & +
  \sum_{X_{i} \in J_{c_{n} \leq \sigma < c_{0}}(\theta)}
  \log{\{f(X_{i};\theta) + \kappa_{0}\}}
  \nonumber \\
  & &   + 
   R_{n}(J_{\sigma < c_{n}}(\theta)) \cdot \log{\frac{v_{0}}{\kappa_{0}}}
  \nonumber \\ 
 & & + 
  \left\{
   R_{n}(J_{\sigma < c_{n}}(\theta)) \cdot \log{\frac{1}{\sigma_{(1)}}}
   + 
   \log{s_{n}(\theta)} 
  \right\}.
  \label{eq:boundingLikelihoodFunc:lem:boundingLikelihoodFunc:thm:penalt-scale-param}
 \end{eqnarray}
\end{lem}

We bound the six terms of (\ref{eq:boundingLikelihoodFunc:lem:boundingLikelihoodFunc:thm:penalt-scale-param}) 
in the following paragraphs. 

\paragraph{The first term}
We begin by bounding the first term of (\ref{eq:boundingLikelihoodFunc:lem:boundingLikelihoodFunc:thm:penalt-scale-param}). 
By lemma~\ref{lem:WaldThm1Type} and the strong law of large
numbers, we have
\begin{eqnarray}
\lim_{n\rightarrow\infty}  \frac{1}{n}\sum_{i=1}^{n}
  \log{
  \left\{
     f_{\scrm_{R}}(X_{i};\theta_{\scrm_{R}},\rho(\theta_{\scrm_{R}})) + 4{\kappa_{0}}
  \right\}
  }
 < 
 E_0[\log f(x;\theta_0)] - 4\lambda_{0} , \quad a.s.
  \label{eq:BoundingTheFirstTerm:thm:penalt-scale-param}
\end{eqnarray}

\paragraph{The second term}
By (\ref{eq:A-zero:condition:thm:penalt-ratio-scale}) and the strong law of large
numbers, we have
\begin{eqnarray}
  \lim_{n\rightarrow\infty} 
  \frac{1}{n}R_{n}(\scra_{0}) \cdot
  \log{
  \left(
   \frac{v_0/c_0 + 3\kappa_{0}}{4\kappa_{0}}
  \right)
  } 
  < 
  {\lambda_{0}}
  , \quad a.s.
  \label{eq:BoundingTheSecondTerm:thm:penalt-scale-param}
\end{eqnarray}

\paragraph{The third term and the fourth term}
The third term and fourth term of (\ref{eq:boundingLikelihoodFunc:lem:boundingLikelihoodFunc:thm:penalt-scale-param}) 
can be bounded from above as follows: 
 \begin{eqnarray}
  & &
   \limsup_{n \rightarrow \infty}
   \sup_{\theta \in \Theta_{c_{n},\scrk,s}^{C}}
   \frac{1}{n} R_{n}(J_{c_{n} \leq \sigma < c_{0}}(\theta))\cdot \abs{\log{2\kappa_{0}}}
   \leq 3M \cdot u_{0} \cdot 2\nu(c_{0}) \cdot \abs{\log{2\kappa_{0}}} < \lambda_{0}, 
   \quad a.s.
   \nonumber \\
   \label{eq:BoundingTheThirdTerm:thm:penalt-scale-param}
   \\
  & &
   \limsup_{n \rightarrow \infty}
   \sup_{\theta \in \Theta_{c_{n},\scrk,s}^{C}}
   \frac{1}{n}  \sum_{X_{i} \in J_{c_{n} \leq \sigma < c_{0}}(\theta)}
   \log{
   \left\{
    f(X_i; \theta) + \kappa_{0}
   \right\}
   }
   \leq
   {\lambda_{0}} ,
   \quad a.s.
   \nonumber \\
  \label{eq:BoundingTheFourthTerm:thm:penalt-scale-param}
 \end{eqnarray}
The proofs of the above inequalities are similar to the proofs of section 4.3.4 and 4.3.5 in longer version of \citet{TT2003-35}, 
and are omitted. 

\paragraph{The fifth term}
We now state the Lemma~\ref{lem:numberOfObservationsWithinShortInterval:thm:penalt-scale-param} 
in order to bound the fifth term of (\ref{eq:boundingLikelihoodFunc:lem:boundingLikelihoodFunc:thm:penalt-scale-param}). 
\begin{lem}
 \label{lem:numberOfObservationsWithinShortInterval:thm:penalt-scale-param}
 \begin{equation}
  \prob\left(
	\sup_{\theta \in \Theta_{c_{n},\scrm,s}^{C}}
	R_{n}(J_{\sigma < c_{n}}(\theta)) > M \quad \io
       \right) = 0 
 \label{eq:lem:numberOfObservationsWithinShortInterval:thm:penalt-scale-param}
 \end{equation}
\end{lem}
\Proof
Let $w_{n} = \nu(c_{n}) = \nu(\exp(n^{-\cnd}))$. 
Then (\ref{eq:assumption:lem:nondense}), the assumption of Lemma~\ref{lem:nondense}, is satisfied. 
From (\ref{eq:def:J:thm:penalt-scale-param}), we have 
\begin{eqnarray}
  \sup_{\theta \in \Theta_{c_{n},\scrm,s}^{C}} 
  R_{n}(J_{\sigma < c_{n}}(\theta)) > M 
  & \quad \Rightarrow \quad & 
  \max_{\mu \in \real}R_{n}([\mu - w_{n},\, \mu + w_{n}] )  > 1 
  \nonumber 
\end{eqnarray}
Therefore, by Lemma~\ref{lem:nondense}, we obtain (\ref{eq:lem:numberOfObservationsWithinShortInterval:thm:penalt-scale-param}). 
\qed

By Lemma~\ref{lem:numberOfObservationsWithinShortInterval:thm:penalt-scale-param} and the same argument in Section~\ref{sec:proof:eq:goal-Phi:thm:penalt-ratio-scale}, 
we ignore the event $R_{n}(J_{\sigma < c_{n}}(\theta)) > M$. 
Then we have
\begin{equation}
 \sup_{\theta \in \Theta_{c_{n},\scrm,s}^{C}}
 R_{n}(J_{\sigma < c_{n}}(\theta)) \cdot \log{\frac{v_{0}}{\kappa_{0}}}
  \leq M \cdot \log{\frac{v_{0}}{\kappa_{0}}}
  \quad a.s. 
  \nonumber 
\end{equation}
Therefore, we obtain for $\theta \in \Theta_{c_{n},\scrm,s}^{C}$
\begin{equation}
 \lim_{n \rightarrow \infty} 
  \sup_{\theta \in \Theta_{c_{n},\scrm,s}^{C}}
  \frac{1}{n} \cdot 
   R_{n}(J_{\sigma < c_{n}}(\theta)) \cdot \log{\frac{v_{0}}{\kappa_{0}}}
   = 
   0
   \quad a.s. 
 \label{eq:BoundingTheFifthTerm:thm:penalt-scale-param}
\end{equation}

\paragraph{The sixth term} 
From Lemma~\ref{lem:numberOfObservationsWithinShortInterval:thm:penalt-scale-param} 
and the same argument in Section~\ref{sec:proof:eq:goal-Phi:thm:penalt-ratio-scale}, 
we have for $\theta \in \Theta_{c_{n},\scrm,s}^{C}$ uniformly 
\begin{equation}
 R_{n}(J_{\sigma < c_{n}}(\theta)) \cdot \log{\frac{1}{\sigma_{(1)}}}
   + 
   \log{s_{n}(\theta)} 
  \leq \log{\frac{s_{n}(\theta)}{(\sigma_{(1)})^{M}}}
  \quad a.s. 
  \nonumber 
\end{equation}
Furthermore, from Assumption~\ref{assumption:10} and \ref{assumption:11}, we have 
\begin{equation}
 \frac{s_{n}(\theta)}{(\sigma_{(1)})^{M}}
  = \frac{\bar{s}_{n}(\sigma_{(1)})}{(\sigma_{(1)})^{M}} \cdot \prod_{m=2}^{M}\bar{s}_{n}(\sigma_{(m)}) 
  \leq \bar{S}^{M-1} \cdot \bar{s} \cdot \exp{(n^{\cnd})}. 
 \nonumber 
\end{equation}
Note that $0 \leq \cnd < 1 $. 
Therefore we obtain for $\theta \in \Theta_{c_{n},\scrm,s}^{C}$
\begin{equation}
 \lim_{n \rightarrow \infty} \frac{1}{n} \cdot 
  \left\{
   R_{n}(J_{\sigma < c_{n}}(\theta)) \cdot \log{\frac{1}{\sigma_{(1)}}}
   + 
   \log{s_{n}(\theta)} 
  \right\}
   = 
   0
   \quad a.s. 
 \label{eq:BoundingTheSixthTerm:thm:penalt-scale-param}
\end{equation}

\paragraph{The end of the proof}
From (\ref{eq:BoundingTheFirstTerm:thm:penalt-scale-param}), 
(\ref{eq:BoundingTheSecondTerm:thm:penalt-scale-param}), 
(\ref{eq:BoundingTheThirdTerm:thm:penalt-scale-param}), 
(\ref{eq:BoundingTheFourthTerm:thm:penalt-scale-param}), 
(\ref{eq:BoundingTheFifthTerm:thm:penalt-scale-param}), 
(\ref{eq:BoundingTheSixthTerm:thm:penalt-scale-param}), 
and Lemma~\ref{lem:boundingLikelihoodFunc:thm:penalt-scale-param}, we have 
\begin{eqnarray}
 \lim_{n \rightarrow \infty}
  \sup_{\theta \in \Theta_{c_{n},\scrm,s}^{C}}\frac{1}{n}
  \left\{
   \sum_{i=1}^{n}\log{f(X_{i};\theta)} + \log{s_{n}(\theta)}
  \right\}
  <
  E_{0}[\log{f(x;\theta_{0})}] - \lambda_{0}
  \quad a.s. 
 \nonumber 
\end{eqnarray}
Therefore we obtain (\ref{eq:goal2:thm:penalt-scale-param}). 

This completes the proof of Theorem~\ref{thm:penalt-scale-param}. 

\section{Conclusion}
\label{sec:conclusion}

In location-scale mixture distributions, 
we have shown the consistency results for 
the two types of penalized maximum likelihood estimators. 
In Corollary~\ref{cor:const-ratio-scale}, an open problem mentioned in \citet{H1985}, \citet{M2000} has been solved positively as follows: 
\begin{itemize}
 \item It is possible to let the lower bound $b$ of the ratios of variances 
       decrease to zero as the sample size $n$ increases to infinity while maintaining consistency. 
 \item If the rate of convergence of $b$ is slower than $\exp(-n^{\bnd})$ where $\bnd$ is a constant such that $0 < \bnd < 1$, 
       then the maximum likelihood estimator is strongly consistent under the constraint $\min_{m, m'} \frac{\sigma_{m}}{\sigma_{m'}} \geq b$. 
\end{itemize}

The assumptions for the penalties given in section~\ref{sec:penal-ratio-scale} or section~\ref{sec:penal-scale-param} are not so restrictive. 
Note that the penalty does not have to depend on the sample size $n$. 
For example, if we set 
$\bar{r}_{n}(y) = \bar{r} \cdot y^{\alpha-1}$ 
and assume $\alpha > M+1$, 
then $\bar{r}_{n}(y)$ satisfies the Assumption~\ref{assumption:6} 
and 
$r_{n}(\theta) = \bar{r}_{n}(\frac{\theta_{(1)}}{\theta_{(M)}})$ satisfies the Assumption~\ref{assumption:8} and \ref{assumption:9}. 
The penalized likelihood $g_{n}(\theta;\vr{X})$ corresponds to the posterior likelihood 
when we adopt a beta distribution as the prior of the minimum of the ratios of the scale parameters. 
Furthermore, if we set $\bar{s}_{n}(y) = {e^{-\frac{\beta}{y}}}\cdot{y^{-(\alpha+1)}}$ and assume $\alpha, \beta > 0$, 
then $\bar{s}_{n}(y)$ satisfies the Assumption~\ref{assumption:10} and \ref{assumption:11}, and 
$s_{n}(\theta) = \prod_{m=1}^{M}\bar{s}_{n}(\sigma_{m})$ satisfies the Assumption~\ref{assumption:12}. 
The penalized likelihood $h_{n}(\theta;\vr{X})$ corresponds to the posterior likelihood 
when we adopt inverse gamma distributions as the priors of the scale parameters. 

From Theorem~\ref{thm:penalt-ratio-scale} and \ref{thm:penalt-scale-param}, 
we can easily show that the consistency of penalized likelihood estimator holds 
even when restrictions on either the location or scale parameters exist. 
If we know that the true parameter is in the restricted parameter space and the assumtions hold, 
then the consistency of the penalized maximum likelihood estimator holds 
by setting $r_{n}(\theta)=0$ or $s_{n}(\theta)=0$ for all $n$ outside the restricted parameter space. 
For example, suppose one considers a uniform mixture distributions under the assumption that the data is non-negative. 
Theorem~\ref{thm:penalt-ratio-scale} and \ref{thm:penalt-scale-param} are applicable to this model. 

\bibliographystyle{econometrica}
\bibliography{mixture}

\end{document}